\documentclass[10pt]{amsart}
\usepackage{amssymb}

\usepackage{graphicx}

\theoremstyle{plain}
\newtheorem{thm}{Theorem}
\newtheorem{prop}[thm]{Proposition}
\newtheorem{lem}[thm]{Lemma}

\newtheorem{cor}[thm]{Corollary}

\newtheorem*{claim}{Claim}
\newtheorem{exm}[thm]{Example}
\newtheorem{rem}[thm]{Remark}
\newtheorem{df}[thm]{Definition}

\renewcommand{\leq}{\leqslant}

\newcommand{\bigjoin}{\bigvee}

\renewcommand{\>}{\rangle}
\newcommand{\rto}{\rightarrow}

\newcommand{\nor}{\triangleleft}

\begin{document}

\title{Ordered direct implicational basis of a finite closure system}
\author {K. Adaricheva}
\address{Department of Mathematical Sciences, Yeshiva University,
New York, NY 10016, USA}
\email{adariche@yu.edu}

\author {J. B. Nation}
\address{Department of Mathematics, University of Hawai'i, Honolulu, HI
96822, USA}
\email{jb@math.hawaii.edu}

\author {R. Rand}
\address{Yeshiva University, New York, NY 10033, USA}
\email{rrand@yu.edu}

\thanks{The first author was partially supported by AWM-NSF Mentor Travel grant N0839954.}

\begin{abstract}
Closure system on a finite set is a unifying concept in logic programming, 
relational data bases and knowledge systems. It can also be presented in 
the terms of finite lattices, and the tools of economic description 
of a finite lattice have long existed in lattice theory. We present this 
approach by describing the so-called $D$-basis and introducing the 
concept of \emph{ordered direct basis} of an implicational system.
A direct basis of a closure operator, or an implicational system, is a 
set of implications that allows one to compute the closure of an arbitrary 
set by a single iteration.
This property is preserved by the $D$-basis at the cost of following a 
prescribed order in which implications will be attended. In particular, 
using an ordered direct basis allows to optimize the \emph{forward chaining procedure} in logic 
programming that uses the Horn fragment of propositional logic.
One can extract the $D$-basis from any direct unit basis $\Sigma$ in time 
polynomial in the size $s(\Sigma)$, and it takes only linear time 
of the cardinality of the $D$-basis to put it into a proper order.
We produce examples of closure systems on a $6$-element set, for which
the canonical basis of Duquenne and Guigues is not ordered direct. 
\end{abstract}

\maketitle
\section{Introduction} 

In K.~Bertet and B.~Monjardet \cite{BM}, it is shown that five implicational bases for a closure operator on a finite set, found in various contexts in the literature, are actually the same.
The goal of this paper is to demonstrate that standard lattice-theoretic results about the ``most economical way" to describe the structure of a finite lattice may be transformed into a basis for a closure system naturally associated with that lattice. 

The coding of a finite lattice in the form of a so-called $OD$-graph was first suggested in \cite{Na}.
We will call the basis directly following from this $OD$-graph a $D$-\emph{basis}, since it is closely associated with a $D$-\emph{relation} on the set of join-irreducibles of a lattice (not necessarily finite) that was crucial in the studies of free and lower bounded lattices, see \cite{FJN}. 
The definition and the proof that $D$-basis does define a given closure system are given in section \ref{DB}. 

The $D$-basis is a subset of a so-called \emph{dependence relation basis} (Definition 6 in \cite{BM}). Thus, it is also a subset of the \emph{canonical direct unit basis} that unifies the five bases discussed in \cite{BM}.
In section \ref{DRB}, we give an example to demonstrate that the reverse inclusion does not hold, thus showing that this newly introduced $D$-basis is generally shorter than the existing ones. 

Recall that the main desirable feature of bases from \cite{BM} is that they be \emph{direct}, which means that the computation of the closure of any subset can be done by attending each implication from the basis only once. This makes the computation of closures a one-iteration process.

While the $D$-basis is not direct in this meaning of this term, 
the closures can still be computed in a single iteration of the basis, 
provided the basis was put in a specific order prior to 
computation. Moreover, there is a simple and effective linear time algorithm 
for ordering a $D$-basis appropriately.  Thus, applying the $D$-basis can be
compared to the iteration known in artificial intelligence as the 
\emph{forward chaining algorithm}, see for example \cite{KL}. 

We introduce the definition of ordered iteration and ordered direct basis in section \ref{ODB}, where
we also prove that the $D$-basis is ordered direct and discuss the algorithmic aspects
of ordering it. 
The further directions of optimization of D-basis are outlined in section \ref{optimizing}, where we also introduce the notion of an \emph{ordered direct sequence} built from a given basis of a closure system.

In section \ref {EB}, we also discuss the so-called $E$-relation, introduced in \cite{FJN}, which leads to the definition of the $E$-basis in closure systems \emph{without $D$-cycles}. In general, the implications written from the $E$-relation do not necessarily form a basis of a closure system, but in closure systems without $D$-cycles, the $E$-basis is ordered direct, is contained in the $D$-basis, and often shorter than the $D$-basis. We discuss a polynomial time algorithm for ordering the $E$-basis.

We explore the connections between $D$-basis, $E$-basis and the so-called \emph{canonical basis} introduced by Duquenne and Guigues in \cite{DG}. While the canonical basis has the minimal number of implications among all the bases of a closure system, it does not have the feature of $D$-basis or $E$-basis discussed in this paper, namely, it cannot be turned into an ordered direct basis. Section \ref{D and DG} of our paper presents examples of closure systems on a $6$-element set, for which the canonical basis cannot be ordered.
As a result, the time required for one iteration of $D$-basis wins over at least two iterations of the canonical basis.
Further polynomial-time optimizations of both $D$-basis and the canonical basis are discussed.

Section \ref{forward} is devoted to discussion and testing the forward chaining algorithm in comparison to the ordered direct basis algorithm. 
Section \ref{performance} provides test results comparing the performance of the $D$-basis with the Duquenne-Guigues canonical basis and canonical direct unit basis. 

The next two sections contain the required definitions and establish connections between finite lattices, closure operators, implicational systems, Horn formulas and Horn Boolean functions.
The reader may consult the survey \cite{CM03} for various aspects of closure systems on finite sets.

\section{Lattices and closure operators}

By a \emph{lattice}, one means an algebra with two binary operations $\wedge, \vee$, called \emph{meet} and \emph{join}, respectively. Both operations are idempotent and symmetric and are connected by absorbtion laws: $x\vee (x\wedge y)=x$ and $x\wedge (x\vee y)=x$. These laws allow us to define a partially order on the base set of the lattice: $x\leq y$ iff $x\wedge y=x$ (which is equivalent to $x\vee y= y$). Vice versa, every partially ordered set, where every two elements have a least upper bound and a greatest lower bound, is, in effect, a lattice. Indeed, in this case the operation $\vee$ can be defined as the least upper bound, and $\wedge$ as the greatest lower bound of two elements in the poset.  A lattice is finite when the base set of this algebra is finite. The symbols $\bigwedge, \bigvee$ are used when more than two elements meet or join.
We will use the notation $0$ for the least element of a lattice, and $1$ for its greatest element.
If $a \leq b$ in lattice $L$, then we denote by $[a,b]$ the interval in $L$, i.e., the set of all $c$ satisfying $a \leq c\leq b$.

Recall now the standard connection between a closure operator on a set and the lattice of its closed sets.
Given a non-empty set $S$ and the set $P(S)=2^S$ of all its subsets, \emph{a closure operator} is a map $\phi: P(S) \rightarrow P(S)$ that satisfies the following, for all $X,Y \in P(S)$: 
\begin{itemize}
\item[(1)] increasing: $X \subseteq \phi(X)$;
\item[(2)] isotone: $X \subseteq Y$ implies $\phi(X)\subseteq \phi(Y)$;
\item[(3)] idempotent: $\phi(\phi(X))=\phi(X)$.
\end{itemize}
It would be convenient for us to refer to the pair $\langle S,\phi\rangle$ of a set $S$ and a closure operator on it as a \emph{closure system}.

A subset $X \subseteq S$ is called \emph{closed} if $\phi(X)=X$. The collection of closed subsets of closure operator $\phi$ on $S$ forms a lattice, which is usually called the \emph{closure lattice} of the closure system $\langle S,\phi\rangle$. This paper deals with only finite closure systems and finite lattices.

Conversely, we can associate with every finite lattice $L$ a particular closure system $\langle S,\phi\rangle$ in such a way that $L$ is isomorphic to a closure lattice of that closure system. Consider $J(L) \subseteq L$, a subset of \emph{join-irreducible elements}. An element $j \in L$ is called join-irreducible, if $j \not = 0$, and $j=a \vee b$ implies $a=j$ or $b=j$. We define a closure system with $S=J(L)$ and the following closure operator:
\[ \phi(X)=[0,\bigvee X]\cap J(L)
\]
It is straightforward to check that the closure lattice of $\phi$ is isomorphic to $L$.

\begin{exm}\label{N5}
Consider a simple example illustrating a closure system built from the lattice $L=\{0,a,b,c,1\}$, for which $0<a<b<1$, $0<c<1$, $a\vee c=b \vee c=1$ and $a\wedge c = b\wedge c=0$. Then $S=J(L)=\{a,b,c\}$. The closed subsets are $[0,x]\cap J(L)$ for $x \in L$, which are $\emptyset$, $\{a\}$, $\{c\}$,$\{a,b\}$ and $\{a,b,c\}$. Knowing all closed subsets, one can define a closure of $X$, or $\phi(X)$, as the smallest closed set containing $X$. For example, $\phi(\{b\})=\{a,b\}$.
\end{exm}

There are infinitely many sets and closure operators whose closure lattice is isomorphic to a given $L$. On the other hand, the one just described is the unique one with two additional properties:
\begin{itemize}
\item[(1)] $\phi(\emptyset)=\emptyset$;
\item[(2)] $\phi(\{i\})\setminus \{i\}$ is closed, for every $i \in S$.
\end{itemize}
Condition (2) just says that each $\phi(\{i\})$ is join irreducible.
Note that (1) is a special case of (2), and that (2) implies the property
\begin{itemize}
\item[(3)] $\phi(\{i\})=\phi(\{j\})$ implies $i=j$, for any $i,j \in S$.
\end{itemize}

Note that $\bigcup_{i \in X} \phi(\{i\}) \subseteq \phi(X)$, but the inverse inclusion does not necessarily hold. In Example \ref{N5}, for instance, $\phi(\{a\})\cup \phi(\{c\})\subset \phi(\{a,c\})$, since $b$ belongs to the right side and not to the left side.

We will call a closure system with properties (1), (2) above a \emph{standard closure system}.
Closure systems with (2) are called $(T\frac{1}{2})$ closure spaces in Wild~\cite{W}.

A closure system satisfying property (3) is said to be \emph{reduced}.
Note that (3) implies $|\phi(\emptyset)| \leq 1$.
Reduced closure systems correspond to a representation of a lattice $L$ as a closure system on a set $S$ with $J(L) \subseteq S \subseteq L$ and $\phi(X)=[0,\bigjoin X] \cap S$.  
A natural example is the set of principal congruences in the congruence lattice of a finite algebra.  Every standard closure system is reduced, and reduced closure systems form a useful intermediate ground between standard and general systems.

It is straightforward to verify that the standard system is characterized by the property that the set $S$ is of the smallest possible size. In other words, one cannot reduce $S$ to define an equivalent closure system. On the other hand, the reduced systems might have excessive elements in $S$.

\begin{exm}
Consider again lattice $L=\{0,a,b,c,1\}$ from Example \ref{N5}. It will represent the closure lattice on $S_1=\{a,b,c,d\}$, where the closed sets are $\emptyset$, $\{a\}$, $\{c\}$,$\{a,b\}$ and $\{a,b,c,d\}$. Thus, in this representation $J(L)\subset S_1$, and property (2) fails: $\phi(\{d\})\setminus \{d\}=\{a,b,c\}$ is not closed. On the other hand, property (3) holds, thus, it is a reduced closure system. Apparently, $S_1$ can be reduce by element $d$, to get an equivalent representation of Example \ref{N5}.
\end{exm}
 
If the closure system $\langle S,\phi\rangle$ is not reduced, one can modify it to produce an equivalent one that is reduced.  Moreover, there is an effective algorithm for doing so.
Thus, for all practical purposes, one can work with a reduced closure system $\langle U,\mu\rangle$ replacing a given one $\langle S,\phi\rangle$. Slightly more effort yields an equivalent standard closure system $\langle V,\nu \rangle$.  The transition is described as follows.

If $\phi(\emptyset)=A \subseteq S$ in  $\langle S,\phi\rangle$, then define $T=S\setminus A$, and redefine a closure operator: $\tau(Y)=\phi(Y)\setminus A$, for all $Y \subseteq T$.
The closure system $\langle T,\tau \rangle$ satisfies property (1).
As (1) is required for a standard closure system, but not for a reduced system, this step may be omitted if only the latter is sought.

Next define an equivalence relation $\approx$ on $T$ by $x \approx y$ if and 
only if $\tau(x)=\tau(y)$.  Then factor out $\approx$, letting $U=T/\approx$
and $\mu(Y)=\tau(Y)/\approx$ for $Y \subseteq U$.  
Alternately, we could define $U$ to be a set of representatives
for $T/\approx$ and $\mu$ to be the restriction of $\tau$.
Either way, one easily checks that $\mu$ is a well-defined closure operator on $U$, and that the closure lattice of 
$\langle U,\mu \rangle$ is isomorphic to that of $\langle S,\phi \rangle$.  
At this point, $\langle U,\mu \rangle$ is reduced.  Moreover, we can recover the original system $\langle S,\phi \rangle$ by expanding the equivalence classes and adding back in $\phi(\emptyset)$.  If desired, we can now continue to produce an equivalent standard closure system.

Let $V = \{ u \in U : \mu(\{ u \}) \setminus \{ u \} \text{ is closed}\}$,
that is, $u \notin  \mu(\mu(\{ u \}) \setminus \{ u \})$, and for 
$Z \subseteq V$ let $\nu(Z) = \mu(Z) \cap V$.  
It is straightforward to verify that $\langle V,\nu \rangle$ is a closure
system satisfying (1) and (2), and that the lattice of closed sets of
$\langle V,\nu \rangle$ is isomorphic to that of $\langle U,\nu \rangle$.  

For the sequel, we will consider primarily reduced closure systems.  Given an arbitrary closure system, not necessarily reduced, the above reduction can be considered as a setup process to allow us to apply the D-basis and related methods.

\section{The bases of closure systems, Horn formulas and Horn Boolean functions}

If $y \in \phi(X)$, then this relation between an element $y \in S$ and a subset $X \subseteq S$ in a closure system can be written in the form of implication: $X \rightarrow y$. Thus, the closure system $\langle S,\phi\rangle$ can be replaced by the set of implications:
\[ 
\Sigma_\phi = \{X \rightarrow y: y \in S, X \subseteq S \text{ and } y \in \phi(X)\}
\]
Conversely, any set of implications $\Sigma$ defines a closure system: the closed sets are exactly subsets
$Y\subseteq S$ that respect the implications from $\Sigma$, i.e., if $X \rightarrow y$ is in $\Sigma$, and $X \subseteq Y$, then $y \in Y$.

It is convenient to define an implication $X\rto y$ as any ordered pair $(X,y)$, $X \subseteq S$, $y \in S$, especially having in mind its interpretation as a propositional formula, see this section two paragraphs below. On the other hand, from the point of view of closure systems, any single implication $ X\rto x$, with $x \in X$, defines a trivial closure system, where all subsets of $S$ are closed. If such implication is present in the set of implications $\Sigma$, then it can be removed without any change to the family of closed sets that $\Sigma$ defines. We will assume throughout the paper that implications $X\rto x$, where $x \in X$, are not included in the set of implications defining closure systems.

Two sets of implications $\Sigma$ and $\Sigma'$ on the same set $S$ are called \emph{equivalent}, if they define the same closure system on $S$. The term \emph{basis} is used for a set of implications $\Sigma'$ satisfying some minimality condition; thus there may be different types of bases. 

Note that, in general, one can consider implications of the form $X \rto Y$, where $Y$ is not necessarily a one-element subset of $S$. Following \cite{BM}, we will call basis $\Sigma$ a \emph{unit implicational basis} if $|Y|=1$ for all implications $X \rto Y$ in $\Sigma$. We will mostly be concerned with unit implicational bases, except for the discussion of the canonical basis of Duquenne-Guigues and its comparison with $D$-basis and $E$-basis. Given any unit basis, we can always collapse the implications with the same premise into one with all conclusions combined into a single set. This will be called an \emph{aggregated} basis.

For a set of implications $\Sigma = \{ X_1 \rto Y_1, \dots, X_m \rto Y_m \}$, define the \emph{size} by $s(\Sigma) = \sum_{j=1}^m (|X_j|+|Y_j|)$.
This is one convenient measure of the complexity of an implicational system.

In general, implications $X\rto y$, where $X \subseteq S$ and $y \in S$, can be treated as the formulas of propositional logic over the set of variables $S$, equivalent to $y \vee \bigvee_{x \in X} \neg x $.  Formulae of this form are also called \emph{definite Horn clauses}. More generally, Horn clauses are disjunctions of negations of several literals and at most one positive literal. The presence of a positive literal makes a Horn clause \emph{definite}.  A \emph{Horn formula} is a conjunction of Horn clauses.

What is called a \emph{model} of a definite Horn clause in logic programming literature corresponds to a closed set of the closure operator defined by this clause. Indeed, by the definition, a model of any formula is simply a tuple $m\in 2^S$ of zeros and ones assigned to literals from $S$, such that the formula is true (=1) on this assignment.
For the definite Horn clause $X \rto y$, $m$ corresponds to a subset $Y$ of $S$ that is closed  for a closure operator on $S$ defined by $X \rto y$. In fact, $m$ is just the characteristic function of $Y$.

There is also a direct correspondence between Horn formulas and Horn Boolean functions: a Boolean function $f:\{0,1\}^n \rto \{0,1\}$ is called a (\emph{pure} or definite) \emph{Horn function}, if it has some CNF representation given by a (definite) Horn formula $\Sigma$.  The dual definition is sometimes
used in the literature, so that a
Horn function is given by some formula in DNF, whose negation is a Horn formula \cite{CH}. Using either definition, one can translate many results on Horn Boolean functions to the language of closure operators, see more details in \cite{BM}.

Consider a set $\Sigma$ of Horn clauses over some finite set of literals $S=\{x_1,\dots,x_n\}$. If some Horn clause $\alpha$ in $\Sigma$  is not definite, i.e., is of the form $\bigvee_{x \in X} \neg x $, $X \subset S$, and it does not use all literals from $S$, then we could define the set of definite clauses $\Sigma_\alpha =\{X \rto y: y \in S \setminus X\}$. It is easy to observe that the set of models of $\Sigma_\alpha$ consists of all models of $\bigvee_{x \in X} \neg x $ and one additional model, which is a tuple of all ones, representing the set $S$ itself. If some clause $\beta \in \Sigma$ is not definite and uses all the literals from $S$, then we define $\Sigma_\beta = \emptyset$ (another possibility, $x_1 \rto x_1$). Again, the set of models $\Sigma_\beta$, i.e., all tuples of zeros and ones, extends the models of $\beta$ by single tuple of all ones. It follows that the set of definite clauses $\Sigma'$, where each non-definite clause $\alpha$ from $\Sigma$ is replaced by a set of clauses $\Sigma_\alpha$, has the set of models that extends the set of models of $\Sigma$ by a single tuple of all ones. This includes the case when $\Sigma$ has no models, i.e., when it is \emph{inconsistent}. 

This observation allows us to reduce the solution of various questions about sets of Horn clauses to sets of \emph{definite} Horn clauses. Thus, it emphasizes the importance of the study of closure operators on $S$.

One of the important questions in logic programming is whether one clause $\phi$ is a \emph{consequence} of the set (or conjunction) of clauses $\Sigma$. Denoted by $\Sigma \models \phi$, this means that every model of $\Sigma$ is also a model of $\phi$. If $\phi$ and formulas in $\Sigma$ are Horn clauses, then, translating this question to the language of closure systems, one reduces it to checking whether every closed set of a closure system defined by $\Sigma$ respects $\phi$. 

\section{The D-basis}\label{DB}

In this section we are going to define a basis that translates to the language of closure systems the defining relations of a finite lattice developed in the lattice theory framework. One can consult \cite{FJN} for the corresponding notion of a minimal cover and $D$-relation used in the theory of free lattices and lower bounded lattices.

Given a \emph{reduced} closure system $\langle S,\phi\rangle$, let us define two auxiliary relations. The first relation is between the subsets of $S$: we write $X\ll Y$, if for every $x \in X$ there is $y \in Y$ satisfying $x \in \phi(y)$. In Example \ref{N5}, for instance, we have $\{a\} \ll \{b\}$, $\{a,c\} \ll \{b,c\}$ and $\{c\} \ll \{a,c\}$.
Note that $X \subseteq Y$ implies $X \ll Y$.
We also write $X\sim_\ll Y$, if $X\ll Y$ and $Y \ll X$. This is true for $X=\{a,b,c\}$ and $Y=\{b,c\}$ in Example \ref{N5}. 

Several observations are easy.

\begin{lem}
The relation $\ll$ is a quasi-order, and thus $\sim_\ll$ is an equivalence relation on $P(S)$.
\end{lem}

We will denote a $\sim_\ll$-equivalence class containing $X$ by $[X]$. Note that for any two members $X,Y \in [X]$, we have $\phi(X)=\phi(Y)$. 
Indeed, $X\ll Y$ implies $X \subseteq \phi(Y)$ and $\phi(X)\subseteq \phi(Y)$. Inverse inclusion follows from $Y\ll X$.

There is a natural order $\leq_c$ on $\sim_\ll$-classes:
$[X] \leq_c [Y]$ if $X \ll Y$.

\begin{lem}\label{poset}
The relation $\leq_c$ is a partial order on the set of $\sim_\ll$-equivalence classes.
\end{lem}

Each class $[X]$ is ordered itself with respect to set containment.

In Example \ref{N5}, we have that $\{a,b,c\} \sim_\ll \{b,c\}$, and no more subsets are $\sim_\ll$-equivalent to $\{a,b,c\}$. Thus, $[\{b,c\}]$ consists of two subsets, and $\{b,c\}\subseteq \{a,b,c\}$ is the minimal (with respect to the order of containment) subset in that equivalence class of $\sim_\ll$. Also $\{a,c\}\ll \{b,c\}$, whence $[\{a,c\}]\leq_c [\{b,c\}]$.

\begin{lem}\label{min}
If $\langle S,\phi \rangle$ is reduced, then
each equivalence class $[X]$ has a unique minimal element with respect to the containment order.
\end{lem}
\begin{proof}
Let us assume that there are two minimal members $X_1$ and $X_2$ in $[X]$.
Without loss of generality we assume that there is $x \in X_1 \setminus X_2$. 
 Since $X_1\ll X_2 \ll X_1$, we have $x \in \phi(x_2)$ and
$x_2 \in \phi(x_1)$, for some $x_1 \in X_1$, $x_2 \in X_2$.  We cannot have $x=x_1$, because, if so, then $\phi(x)=\phi(x_2)$, which implies $x=x_2$, since our closure system is reduced. 
This would contradict to the choice of $x$.

Thus, $x \ne x_1$, and $x \in \phi(x_1)$. But then we can reduce $X_1$ to $X'=X_1 \setminus x \subset X_1$, which is still a member of $[X]$ since $X_1 \ll X' \subset X_1$.  This contradicts the minimality of $X_1$ in $[X]$.
\end{proof}

The second relation we want to introduce in this section is between an element $x \in S$ and a subset $X \subseteq S$, which will be called a \emph{cover} of $x$.
(In lattice theory, the terminology \emph{nontrivial join cover} is used.)
We will write $x \nor X$, if $x \in \phi(X)\setminus \bigcup_{x' \in X}\phi(x')$.
This notion is illustrated in Example~\ref{N5} by $b\nor \{a,c\}$. Note that it is not true that $a\nor \{b,c\}$, because $a <b$, so that $a\in \phi(b)$ for the corresponding standard closure operator.

We will call a subset $Y\subseteq S$ a \emph{minimal cover} of an element $x \in S$, if $Y$ is a cover of $x$, and for every other cover $Z$ of $x$, $Z \ll Y$ implies $Y \subseteq Z$.
So a minimal cover of $x$ is a cover $Y$ that is minimal with respect to
the quasi-order $\ll$, and minimal with respect to set containment within its 
$\sim_\ll$-equivalence class $[Y]$, as per Lemma~\ref{min}.

To illustrate this notion, let us slightly modify Example \ref{N5}. Rename element $0$ by $d$ and add a new $0$ element: $0<d$, resulting in a lattice $L_1$ with $J(L_1)=J(L)\cup \{d\}$. We will have $Y=\{a,c\}$ as a minimal cover for $b$. Indeed, the only other cover for $b$ is $Z=\{a,c,d\}$, for which we have $Z\ll Y$ and $Y\subseteq Z$. 

\begin{lem}\label{refines}
For a reduced closure system, if $x \nor X$, then there exists $Y$ such that $x \nor Y$, $Y \ll X$ and $Y$ is a minimal cover for $x$.
In other words, every cover can be $\ll$-reduced to a minimal cover.
\end{lem}
\begin{proof}
Consider $P_x=\{[X]: x\nor X\}$, a sub-poset in the $\leq_c$ poset of $\sim_\ll$-classes.
If it is not empty, choose a minimal element in this sub-poset, say $[Y]$, and let $Y$ be the unique minimal element in $[Y]$ with respect to containment, which exists due to Lemma \ref{min}. Then $Y \ll X$ and $x \nor Y$.
It remains to show that, for every other cover $Z$ of $x$, $Z \ll Y$ implies $Y \subseteq Z$.
Indeed, since $Z \ll Y$, we have $[Z]\leq_c [Y]$. But $[Y]$ is the minimal element in $P_x$, hence, $[Z]=[Y]$.
It follows that $Y \subseteq Z$, since $Y$ is the minimal element of $[Y]$ with respect to containment order.
\end{proof}

We finish this section by introducing the $D$-basis of a reduced closure system.

\begin{df}
Given a reduced closure system $\langle S,\phi\rangle$, we define the $D$-basis $\Sigma_D$ as a union of two subsets of implications:
\begin{itemize}
\item[(1)] $ \{y \rightarrow x : x \in \phi(y)\setminus y,\  y \in S\}$;
\item[(2)] $\{X \rightarrow x: X \text{ is a minimal cover for } x\}$.
\end{itemize}
\end{df}

Part (1) in the definition of the $D$-basis will also be called the \emph{binary part} of the basis, due to the fact that both the premise and the conclusion of implications in (1) are one-element subsets of $S$.

For the closure system $\langle J(L), \phi\rangle$ associated with the lattice $L$ in Example \ref{N5}, the $D$-basis consists of two implications: $b \rightarrow a$ and $\{a,c\} \rightarrow b$. 

\begin{lem}
$\Sigma_D$ generates $\langle S,\phi\rangle$.
\end{lem}
\begin{proof}
We need to show that, for any $x \in S$ and $X \subseteq S$ such that $x \in \phi(X)\setminus X$, the implication $X \rightarrow x$ follows from implications in $\Sigma_D$.


If $x \in \phi(\{x'\})$, for some $x' \in X$, $x'\not = x$, then  $X \rightarrow x$ follows from $x'\rightarrow x$ that is in $\Sigma_D$.
So assume that $x \not \in  \phi(\{x'\})$, for any $x'\in X$. Then $x \nor X$. According to Lemma~\ref{refines}, there exists $Y\ll X$ such that $x \nor Y$, and $Y$ is a minimal cover for $x$. Then $Y \rightarrow x$ is in $\Sigma_D$.
Besides, for each $y \in Y\setminus X$ there exists $x_y \in X$ such that $y \in \phi(\{x_y\})$.  Therefore, $x_y \rightarrow y$ is in $\Sigma_D$ as well. Evidently, $X \rightarrow x$ is a consequence of $Y \rightarrow x$ and $\{x_y\rightarrow y: y \in Y\}$.
\end{proof} 

\section{Comparison of the $D$-basis and the dependence relation basis}
\label{DRB}

One of the bases discussed in \cite{BM} is the \emph{dependence relation basis}. 
For a closure system $\langle S,\phi \rangle$, not necessarily reduced, the dependence relation basis is 
\[
\Sigma_\delta=\{X \rightarrow y: y \in \phi(X) \setminus X \text{ and } y \notin \phi(Z) \text{ for all } Z \subset X \}.
\]
Since $Z \subseteq X$ implies $Z \ll X$, a minimal cover (as defined above) is automatically minimal with respect to containment. Thus we have the following connection.

\begin{lem}\label{subset}
For a reduced closure system, $\Sigma_D \subseteq \Sigma_\delta$.
\end{lem}

For later reference, the dependence relation $\delta$ from Monjardet~\cite{M} can be described by $y \delta x$ whenever $x \in X$ for some $X \rto y$ in $\Sigma_\delta$.  

In the next example and in the sequel, whenever there is no confusion, we will omit the braces in notations of subsets of some set $S$: $\{x\}, \{a,b,c\}$, etc.~will be denoted simply $x$, $abc$, etc.

\begin{figure}[htbp]
\begin{center}
\includegraphics[height=1.7in,width=6.0in]{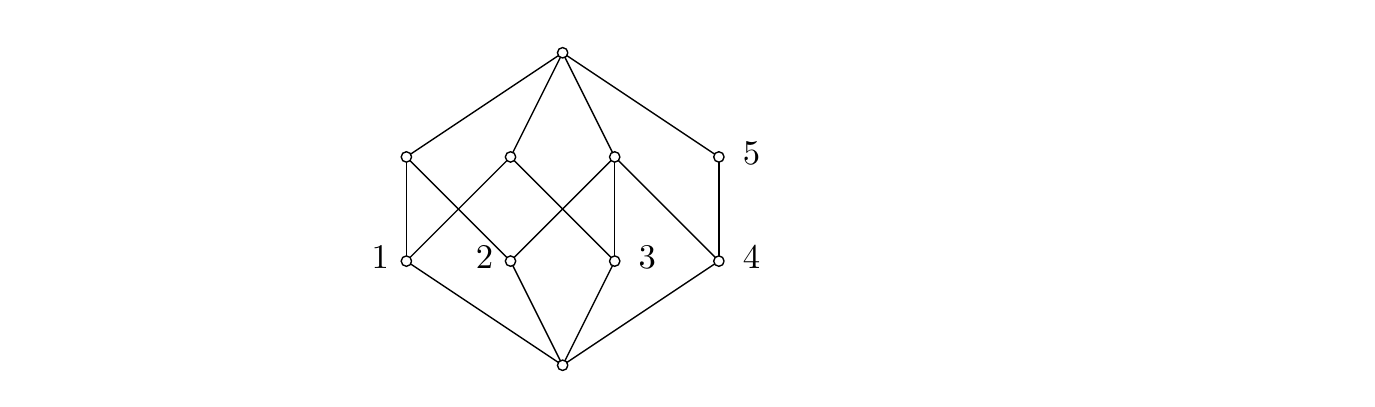}
\caption{Example \ref{ex9}}
\end{center}
\end{figure}

\begin{exm} \label{ex9}
This example is based on Example 5 from \cite{BM}.
Consider the closure system on $S=\{1,2,3,4,5\}$ with the set of closed subsets $F=\{\emptyset,1,2,3,4,12,13,$ $234,45,12345\}$. Then $\Sigma_\delta=\{5\rightarrow 4, 23\rightarrow 4, 24\rightarrow 3,$ $34\rightarrow 2, 14\rightarrow 2, 14\rightarrow 3, 14\rightarrow 5, 25\rightarrow 1, 35\rightarrow 1,$ $15\rightarrow 2, 35\rightarrow 2, 15\rightarrow 3, 25\rightarrow 3, 123\rightarrow 5\}$.
\end{exm}

All implications except $5\rightarrow 4$ are of the form $X\rightarrow x$, where $x\nor X$. On the other hand, not all covers $X$ are minimal covers of $x$. We can check that each of implications $15\rightarrow 2, 35\rightarrow 2, 15\rightarrow 3, 25\rightarrow 3$ does not represent a minimal cover. For example, $2 \nor 15$, but $14 \ll 15$ and $2\nor 14$ is the minimal cover. In particular, $D$-basis consists of all implications from $\Sigma_\delta$ except the four indicated: $\Sigma_D=\{5\rightarrow 4, 23\rightarrow 4, 24\rightarrow 3, 34\rightarrow 2, 14\rightarrow 2, 14\rightarrow 3, 14\rightarrow 5, 25\rightarrow 1, 35\rightarrow 1, 123\rightarrow 5\}$. 

As this example demonstrates, the $D$-basis can be obtained from $\Sigma_\delta$ simply by removing some unnecessary implications. It turns out that the same can be done for the big range of bases called \emph{direct unit bases}. Moreover, it can be done in polynomial time in the size of the given basis. See Proposition \ref{extractD} in the next section.

\section{Direct basis versus ordered direct basis}\label{ODB}

The bases discussed in Bertet and Monjardet \cite{BM} are, in general, redundant: a proper subset of such a basis would generate the same closure system. For example, as we saw in the previous section, $\Sigma_\delta$ from Example 5 was reduced to a smaller basis $\Sigma_D$.  Example~\ref{ex66}
shows that the D-basis can also be redundant; see Remark~\ref{rem29}. 

While the desire to keep the basis as small as possible might be a plausible task, there is another property of a basis that could be better appreciated in a programming setting. Here we recall the definition of a \emph{direct basis}.

If $\Sigma$ is some set of implications, then let $\pi_\Sigma(X)=X \cup \bigcup \{B: A \subseteq X \text{ and } (A\rightarrow B)\in \Sigma\}$. In order to obtain $\phi_\Sigma(X)$, for any $X \subseteq S$, one would normally need to repeat several iterations of $\pi$: $\phi(X)=\pi(X)\cup \pi^2(X) \cup \pi^3(X) \dots $. 

The bases for which one can obtain the closure of any set $X$ performing only one iteration, i.e., $\phi(X)=\pi(X)$, are called \emph{direct}.

It follows from Theorem 15 of \cite{BM} that the dependency relation basis $\Sigma_\delta$ is direct. Moreover, this basis is direct-optimal, meaning that no other direct basis for the same closure system can be found of smaller total size.
(The \emph{total size} $t(\Sigma)$ is the sum of the cardinalities of all sets participating in its implications.  This will be less than $s(\Sigma)$ if some sets are repeated.)  In particular, any reduction of $\Sigma_\delta$ will cease to be direct. Thus, there is a apparent trade-off between the number of implications in the basis and the number of iterations one needs to compute the closures of subsets.

The goal of this section to implement a different approach to the concept of iteration. That would allow the same number of programming steps as with the iteration of $\pi$, while allowing us to reduce the bases to a smaller size.

\begin{df}\label{ordered iteration}
Suppose the set of implications $\Sigma$ is equipped with some linear order $<$, or equivalently, the implications are indexed as $\Sigma=\{s_1,s_2,\dots, s_n\}$.
Define a mapping $\rho_\Sigma: P(S) \rightarrow P(S)$ associated with this ordering as follows. For any $X \subseteq S$, let $X_0=X$. If $X_k$ is computed and implication $s_{k+1}$ is $A \rightarrow B$, then 
\[
X_{k+1}=\left\{\begin{array}{ll}X_k \cup B, & {\rm if} \; \; A \subseteq X_k,
\\ X_k, & {\rm otherwise.}\end{array}\right.
\]
Finally, $\rho_\Sigma(X)=X_n$.
We will call $\rho_\Sigma$ an \emph{ordered iteration} of $\Sigma$.
\end{df}

\vspace{0.3cm}
Apparently, $\pi_\Sigma(X) \subseteq \rho_\Sigma(X)$, because all implications from $\Sigma$ are applied to the original subset $X$, while they are applied to potentially bigger subsets $X_k$ in the construction for $\rho_\Sigma(X)$.
We note though that assuming the order on $\Sigma$ is established, the number of computational steps to produce $\rho_\Sigma(X)$ is the same as for $\pi_\Sigma(X)$.

\begin{df} The set of implications with some linear ordering on it, $\langle \Sigma, <\rangle$, is called an \emph{ordered direct basis}, if, with respect to this ordering, $\phi_\Sigma(X)=\rho_\Sigma(X)$ for all $X \subseteq S$.
\end{df}

Our next goal is to demonstrate that $\Sigma_D$ is, in fact, an  ordered direct basis. Moreover, it does not take much computational effort to impose a proper ordering on $\Sigma_D$.

\begin{thm}\label{ordering}
Let $\Sigma_D$ be the D-basis for a reduced closure system.
Let $<$ be any linear ordering on $\Sigma_D$ such that all implications of the form $y \rightarrow z$ precede all implications of the form $X\rightarrow x$, where $X$ is a minimal cover of $x$. Then, with respect to this ordering, $\Sigma_D$ is an ordered direct basis. 
\end{thm}

\begin{proof}
Suppose that $X \subseteq S$ and $b \in \phi(X)\setminus X$.  We want to show that $b$ will appear 
in one of the $X_k$ in the sequence that leads to $\rho(X)$.

If $b \in \phi(\{a\}) \setminus \{a\}$ for some $a \in X$, then $b$ will appear in some $X_k$, when $a\rightarrow b$ from $\Sigma_D$ is applied. So now assume that $b \notin \phi(\{a\})$ for every $a \in X$. Then $b \nor X$ and, according to Lemma~\ref{refines}, there exists $Y\ll X$ such that $b \nor Y$ and $Y$ is a minimal cover for $y$. It follows that for any $y \in Y$ there exists $a \in X$ such that $y \in \phi(a)$. All implications $a \rightarrow y$ will be applied prior to any application with the minimal cover. It follows that by the time the implication $s_{k}$, say $Y \rightarrow b$, is tested against $X_{k-1}$, we will have $Y \subseteq X_{k-1}$. Hence, $X_{k}=X_{k-1}\cup\{b\}$. 
\end{proof}

\begin{cor}
$D$-basis is also ordered direct in its aggregated form.
\end{cor}
Indeed, it follows from the fact that the only restriction on the order of the $D$-basis is to have its binary part prior to the rest of the basis. 

\begin{cor}
If\/ $\Sigma_D=\{s_1,\dots,s_m\}$ is the $D$-basis of a reduced implicational system $\Sigma$, then it requires time $O(m)$ to turn it into an ordered direct basis of $\Sigma$. 
\end{cor}

\begin{exm}\label{Dbas}
Consider the closure system with $S=\{1,2,3,4,5,6\}$ and the family of closed sets $F=\{1,12,13,4,45,134,136,1362,1346,13456,123456\}$. Then the $D$-basis of this system is
$\Sigma_D=\{5\rto 4, 14\rto 3, 23\rto 6, 6\rto 3, 15\rto 6,24 \rto 6, 24\rto 5, 3\rto 1, 2\rto 1\}$.
According to Theorem \ref{ordering}, a proper ordering that turns this basis into ordered direct can be defined, for example, as: (1) $5\rto 4$, (2) $6\rto 3$, (3)  $3\rto 1$, (4)  $2\rto 1$, (5) $14\rto 3$, (6) $23\rto 6$, (7) $15\rto 6$, (8) $24 \rto 6$, (9) $24\rto 5$.
\end{exm}

\begin{figure}[htbp]
\begin{center}
\includegraphics[height=2.8in,width=6.0in]{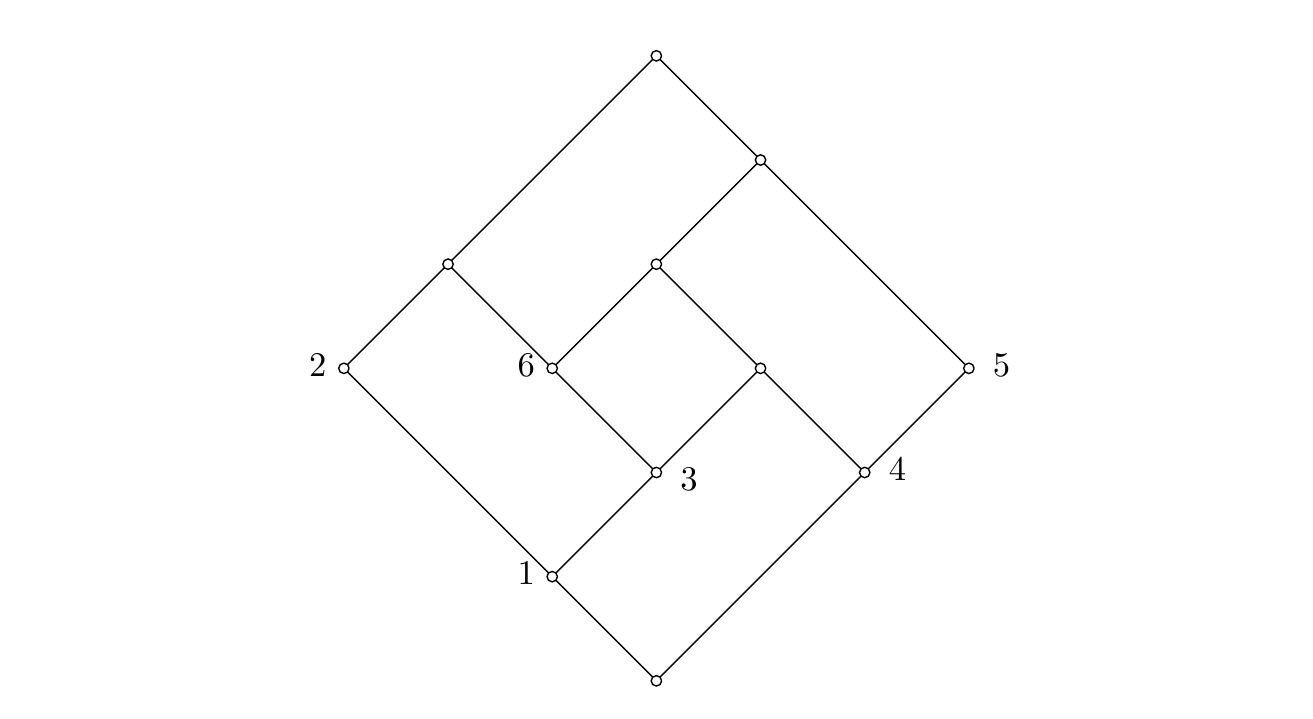}
\caption{Example \ref{Dbas}}
\end{center}
\end{figure}



\section{Processing of ordered basis versus forward chaining algorithm}\label{forward}

The \emph{forward chaining algorithm} was originally introduced in 1984 by W.~Dowling and J.H.~Gallier in the context of checking the satisfiability of Horn formulae \cite{DG84}. In 1992, H.~Mannila and K.J.~R\"{a}ih\"{a} \cite{MR} introduced the \emph{LINCLOSURE} algorithm, which applies the same approach to expanding functional dependencies in Database Systems. In this section, we will look at the efficiency of this approach in comparison with a \emph{folklore} algorithm to computing closures for the $D$-Basis, an approach that can be generalized to any direct or ordered basis.  

We will assume that the base set is $S=\{x_1,\dots, x_n\}$, which can be interpreted as propositional variables, and the closure system is given by a unit basis $\Sigma =\{A_1\rightarrow b_1,\dots, A_m\rightarrow b_m\}$.

The forward chaining procedure requires a pre-processing setup, during which it constructs three data structures:
\emph{ClauseList}$_i=\{A_j : x_i \in A_j\}$ for $i\leq n$,
\emph{Propositions}$_j = |A_j|$ and
\emph{Consequent}$_j = \{b_j\}$ for $j\leq m$,
along with subset \emph{True} $\subseteq S$ thought of as an input set, whose closure needs to be computed. 

When forward chaining computes the closure of \emph{True}, for each new $x_i \in$ \emph{True},  for each $A_j \in$ \emph{ClauseList}$_i$, it decrements the value of \emph{Propositions}$_j$ by one. Whenever \emph{Propositions}$_j=0$, \emph{Consequent}$_j$ is added to the set \emph{True}.

Since every entry of \emph{Propositions} will, in the worst case, be reduced to zero, the number of steps in computing the closure is bounded by the size $s(\Sigma)$, i.e., the combined length of the implications in the basis.  
Including the pre-processing steps, the forward chaining algorithm should require $O(s(\Sigma))$ operations to compute the closure. 
If the closures of multiple sets are to be performed, of course, the setup steps can be abbreviated:  only \emph{Propositions} and \emph{True} need to be updated for subsequent runs.

As noted in \cite{W95}, forward chaining, while efficient in the worst case, generally underperforms the \emph{folklore} algorithm of simply checking if each $A_j$ is contained within \emph{True}, and if so appending $b_j$ to \emph{True}, until the ability of the algorithm to generate new \emph{True} elements is exhausted. In particular, forward chaining does poorly on large sets, where we often only need to examine a fraction of the $|\Sigma|$ variables examined in the forward chaining procedure.

As an alternative to the forward chaining procedure, M.~Wild \cite{W95} suggested an algorithm that considers the set difference $\Sigma'=\Sigma \setminus \{A_k\rightarrow b_k : A_k \not \subseteq \text{ \emph{True}}\}$. For each $(A_j \rto b_j)\in \Sigma'$ it then adds $b_j$ to \emph{True} and repeats as necessary. This algorithm retains the need for preprocessing in the form of \emph{ClauseLists}. Though typically faster than forward chaining, Wild's algorithm has a worst-case running time of $O(s(\Sigma)m^2)$, which can cause problems for large values of $m$.

Applying the folklore algorithm to processing ordered bases, theoretically, avoids the pitfalls of both forward chaining and Wild's algorithm. It simply iterates from $(A_1 \rto b_1)$ to $(A_m \rto b_m)$ adding $b_i$ to \emph{True} whenever $A_i \subseteq \emph{True}$. On one hand, its worst case processing time is $O(s(\Sigma))$ since we only need to iterate through the ordered basis once. At the same time, it takes a fraction of the time of the folklore approach on an non-ordered basis, which will require a minimum of two iterations in order to confirm that no new variables were added to \emph{True}. 

In testing the performance of these three algorithms, we generated $D$-bases from the domains $\{1,2,3,4,5\}$ through $\{1,2,3,4,5,6,7,8\}$ and calculated the time necessary to derive the closure of some random subset of the set. For forward chaining and Wild's algorithm, which require substantial preprocessing, we calculated the time with and without the preprocessing, which corresponds to the time required for computing the first closure on a basis versus the time for subsequent closures.

In our testing, the folklore algorithm considerably outperformed both forward chaining and Wild's algorithm (without preprocessing) though its advantage fell, as the domain grew larger. For a domain of size 5, forward chaining and Wild's algorithm took 2 and 2.66 times as long as folklore, respectively, which shrank to 1.45 and 1.76 as the domain grew to size 8. If we include preprocessing times, however, both algorithms continued to take over 4 times as long, with the relative time of forward chaining remaining relatively constant. Domains of sizes 5-8 corresponded to bases with an average of 8, 13, 19 and 27 implication, respectively.

Taking into an account that LINCLOSURE and Wild's algorithm are normally performed on the agregated bases, we also ran a series of similar tests with the aggregated bases. Such a test is based on the fact that the $D$-basis is ordered direct in both the unit and the aggreagted form. The bases on domains of sizes 5-8 had an overage of 5,7,10 and 13 implications, respectively. The test showed even higher ratios of forward chaining (without preprocessing) times to the ordered direct processing times: from  3.04 for domain 5 to 2.1 for domain 8.

\begin{table}
\begin{tabular}{|l|l|l|l|l|}
\hline
Average Time ($\mu s$) & Domain 5 & Domain 6 & Domain 7 & Domain 8 \\ \hline
Folklore & 3.87 & 6.59 & 10.12 & 14.90 \\ \hline 
Forward Chaining (preprocessed) & 7.74 & 11.32 & 15.99 & 21.58 \\  \hline 
Forward Chaining & 30.34 & 46.10 & 66.74 & 91.48 \\  \hline 
Wild's Algorithm (preprocessed) & 10.14 & 14.51 & 19.66 & 26.16 \\  \hline 
Wild's Algorithm & 23.75 & 35.28 & 50.11 & 69.58 \\ \hline 
\end{tabular}
\caption{Comparing Algorithm Processing Times}
\end{table}



Noticeably, the ordered-basis approach does not actually require the representation of propositions as $S=\{x_1,\dots, x_n\}$ and implications as $\Sigma=\{s_1,\dots, s_m\}$, where each proposition has an associated integer value, necessary for indexing and traversing \emph{ClauseList}, \emph{Propositions}, and \emph{Consequent}. Though we can in principle take advantage of integer values in constructing our set of true values, we only require \emph{a set} of satisfied propositions. By contrast, to use the forward chaining method on a basis without this representation would require significant overhead in hashing each proposition to its corresponding integer.

Additionally, the ordered-basis approach eliminates the need for pre-processing of the basis to store it in the form of \emph{ClauseList} and
\emph{Consequent}. Since the $D$-basis is defined as the union of binary and non-binary sets of implications, which is reflected in the algorithm for producing it, we assume all $D$-bases are properly ordered. This is particularly important when the basis may not fit into main memory. Instead of having to individually access each  \emph{ClauseList}$_i$ when the propositional variable $x_i$ appears in \emph{True}, the ordered-basis approach allows us to parse the basis in conveniently sized pieces.

There is at least one observation how the idea of the ordered basis may improve the performance of the forward chaining algorithm. Indexing the implications according to the proper order of the $D$-basis, whenever we add a variable to \emph{True}, we may additionally maintain the index $j$ of the implication from which it was derived. Then, when we process this variable, we only need to update $k$-entries of \emph{Propositions} where $k\geq j$, saving us significant processing time for very large sets.

\section{Building and optimizing the $D$-basis}\label{optimizing}

We consider the $D$-basis a good alternative of any direct basis, since it has a smaller size than any direct basis and preserves the directness property, under a special ordering we define. In this section we consider an effective procedure to obtain the $D$-basis from any given direct basis, also to further optimize its binary part or to use the concept of the ordered sequence. The problem of obtaining the $D$-basis from other \emph{non-direct} bases is also tackled in \cite{JK2}.

As we saw in Lemma \ref{subset} and Example \ref{ex9}, the $D$-basis $\Sigma_D$ of any reduced closure system is a subset of the direct unit basis $\Sigma_\delta$. The next statement shows that, given any direct unit basis, one can extract the $D$-basis from it in a polynomial time procedure.
 
\begin{prop}\label{extractD}
Let $\langle S,\phi \rangle$ be a reduced closure system.
If the direct unit basis $\Sigma$ for this system has  $m$ implications, and $|S|=n$, then it requires time $O((nm)^2)\sim O(s(\Sigma)^2)$ to build the $D$-basis $\Sigma_D$ equivalent to $\Sigma$.
\end{prop}
\begin{proof}
Let $\langle S,\phi \rangle$ be the closure system on set $S$ defined by $\Sigma$.
By Lemma~\ref{subset}, $\Sigma_D\subseteq \Sigma_\delta$. According to Theorem~15 of \cite{BM}, $\Sigma_\delta$ coincides with the canonical iteration-free basis introduced by M. Wild in \cite{W}. Hence, by Corollary~17 of \cite{BM}, $\Sigma_\delta$ is the smallest basis, with respect to containment, of all direct unit bases of $\langle S,\phi \rangle$. Therefore, $\Sigma_D \subseteq \Sigma_\delta \subseteq \Sigma$. 

It follows that $\Sigma_D$ can simply be extracted from $\Sigma$ by removing unnecessary implications.
This amounts to finding the implications $X \rto x$, where $X$ will be a minimal join cover of $x$, among the implications of $\Sigma$. 

Note that $O(m)$ steps will be needed to separate binary implications $y\rto x$ from $X\rto x$, where $|X|>1$. The number of $x\in S$ that appear in the consequence of implications $X \rto x$ is at most the minimum of $m$ and $n$. 

For every fixed $x$, it will take time $O(m)$ to separate all implications $X\rto x$, and the number of such implications is at most $m$. If $X_1\rto x$ and $X_2\rto x$ are two implications in this set, we can decide in time $O(mn)$ whether $X_1 \ll X_2$ or $x \in \phi(y)$ for some $y \in X_2$. If either holds, $X_2 \rto x$ does not belong to the $D$-basis.  To check this, consider the closure systems $\Sigma_i\subseteq \Sigma$, $i=1,2$ that consist of all binary implications of $\Sigma$, in addition to $X_i\rto x$. Also, put an order on $\Sigma_i$, where all the binary implications precede $X_i\rto x$.
Apparently, $x$ is in the closure of $X_2$, in the closure system defined on $S$ by $\Sigma_1$, iff either $X_1 \ll X_2$ or $x \leq y$ for some $y \in X_2$.

As pointed out in section \ref{forward}, computation of the closure of any input set, either by the forward chaining algorithm, or by the ordered basis algorithm, is linear in the size of the input, which in this case is essentially the size of the binary part of $\Sigma$, or $O(n^2)$.

At the worst case, about $O(m^2)$ comparisons have to be made, for different covers $X_1,X_2$ of the same element $x$, to determine the minimal ones. Hence, the overall complexity is $O(m^2n^2)\sim O(s(\Sigma)^2)$.
\end{proof}

It follows from the procedure of Proposition \ref{extractD} that the $D$-basis is obtained from any direct unit basis by removing implications $X\rto x$, for which $X$ is not a minimal cover of $x$ and $|X|>1$. In particular, the binary part of the direct basis, i.e., implications of the form $y\rto x$, remain in the $D$-basis.

We want to discuss a further optimization of the $D$-basis, as well as any other basis that has the same binary part as the $D$-basis.  As was observed in section~2, for a reduced closure system $\langle S,\phi \rangle$, the elements of $S$ can be identified with elements of the closure lattice $L$, in such a way that $J(L) \subseteq S \subseteq L$.     
This correspondence induces a natural order on $S$, with $s \leq t$ if and
only if $\phi(s) \subseteq \phi(t)$.
Thus, an implication $y\rto x$ belongs to the $D$-basis iff $x \in \phi(y)$ iff $x \leq y$.   The binary part of the $D$-basis then describes the partially ordered set $(S,\leq)$.

Recall that, in the language of ordered sets, we say that $y$ \emph{covers} $x$ if $y>x$ and there is no element $z$ such that $y>z>x$.

We can shorten the binary part of the $D$-basis, leaving only those implications $y\rto x$ for which $y$ covers $x$ in $(S,\leq)$. This will come at the cost of the need to order the remaining implications.
For example, if $x\rto y$, $y\rto z$, $x\rto z$ are three implications from the binary part of some $D$-basis, then the last implication can be removed, under condition that the first two will be placed in that particular order into the ordered $D$-basis. 
More generally, suppose only covering pairs are to be included in the binary part of an ordered basis. Then the ordering of the implications should be such that, if $x>y\geq z>t$ in $S$ with the strict inequalities being covers, then $x \rto y$ precedes $z \rto t$.

Recall also that if some set of implications $\Sigma'$ is ordered, then $\rho_{\Sigma'}(X)$, the ordered iteration of $\Sigma'$, is defined for every $X \subseteq S$, see Definition~\ref{ordered iteration}.

\begin{prop}\label{cover}
Let $\Sigma_1$ be the binary part of the $D$-basis of a reduced closure system on a set $S$.  If\/ $\Sigma_1$ has $k$ implications and $|S|=n$, then there is an $O(nk+n^2)$ time algorithm that extracts $\Sigma'\subseteq \Sigma_1$ describing the cover relation of join irreducible elements of closure system, and places the implications of\/ $\Sigma'$ into a proper order. Under this order, $\rho_{\Sigma'}(y)=\rho_{\Sigma_1}(y)$ for every $y \in S$.
\end{prop}

\begin{proof}
We have the partially ordered set $(S,\leq)$ of size $n$, whose cover relation has at most $k$ pairs,
 thus,  it will take time $O(nk+n^2)$ to find the cover relation of this poset, see \cite{GK79}, also Theorem 11.3 in \cite{FJN}. Let $\Sigma'\subseteq \Sigma_1$ be the set of all implications $y\rto x$, where $y$ covers $x$ in $(S,\leq)$ . It remains to put these implications into a proper order. 
If $(S,\leq_1)$ is any linear extension of $(S, \leq)$, then one can take any order of $\Sigma'$ associated with this extension.  Starting from the maximal element $y$ of  $(S, \leq_1)$, write all implications $y\rto x$ from $\Sigma'$, in any order, then pick next to maximal element $z$ of $(S, \leq_1)$ and write all implications $z \rto t$, in any order, then proceed with all elements of $(S, \leq_1)$ in the same manner, in descending order $\geq _1$. It remains to notice that there is an $O(n+k)$ algorithm for producing the linear extension of partially ordered set with $n$ elements and $k$ pairs of comparable elements, see Theorem 11.1 in \cite{FJN}. 
\end{proof}

Now we want to deviate slightly from the notion of ordered direct basis to the notion of \emph{ordered direct sequence of implications}. Suppose $\Sigma$ is some basis of a closure system $\langle S,\phi\rangle$. The ordered sequence $\sigma=\langle s_1,\dots, s_t\rangle$ of implications from $\Sigma$, not all necessarily different, is called \emph{an ordered direct sequence from }$\Sigma$, if $\rho_\sigma (X)=\phi(X)$ for every $X\subseteq S$.
 
The idea of ordered direct sequencing allows some further optimization of the $D$-basis.
If $Z=\langle z_1,\dots, z_k\rangle$ and $T=\langle t_1,\dots, t_s\rangle$ are two ordered sequences, then $Z^\frown T$ denotes their concatenation (the attachment of $T$ at the end of $Z$). 

\begin{lem}\label{replacing}
Suppose $\sigma=\Sigma_1^\frown\Sigma_2^\frown\Sigma_3$ is an ordered direct sequence from some basis $\Sigma$, where $\Sigma_1$, $\Sigma_3$ consist of binary implications in proper order of Proposition \ref{cover}, $\Sigma_2$ consists of non-binary implications, and $\Sigma_2$ can be put into arbitrary order without changing the ordered direct status. If $(A\rto y),(A\rto x) \in \Sigma_2$ and $(y\rto x) \in \Sigma_1$, then $A\rto x$ can be dropped from $\Sigma_2$ and replaced by an additional $y\rto x$ in $\Sigma_3$.
\end{lem}

\begin{proof}
We need to show that whenever $Y$ is an input set such that $x\in \phi(Y)$, the replacement of $A\rto x$ by $y\rto x$ will not affect computation of $\rho_\sigma(Y)$. 

Consider the case when $y\not \in \phi(Y)$. Then also $A\not \subseteq \phi(Y)$, whence any implication with the premise $A$ will never be applied in computation of $\rho_\sigma(Y)$. The same is true for implications with premise $y$, so replacement of $A\rto x$ by $y\rto x$ can trivially be done.

Now suppose that $y \in \phi(Y)$.  By assumption, we can take $A\rto x$ to be the last implication in the ordering of $\Sigma_2$. So consider $Y_k$, the result of ordered iteration of $\Sigma_1^\frown\Sigma_2\setminus (A\rto x)$ on the input set $Y$. If $y \in Y_k$, then we can drop $A\rto x$ from $\Sigma_2$ and place $y\rto x$ anywhere in proper order in $\Sigma_3$, which will guarantee that $x$ appears in $\rho(Y)$.
If $y \not \in Y_k$, then there is $z \in Y_k$ such that there exists some sequence in $\Sigma_3$ from $z$ to $y$. 
By assumption, $\Sigma_3$ is in the proper order, hence any implication $w\rto y$ precedes $x\rto t$. Thus, we can place $y\rto x$ in between those groups, following the proper order on all binary implications from Proposition \ref{cover}. After replacing $A\rto x$ by $y\rto x$ in proper position of $\Sigma_3$, we can still assume that the ordering of remaining part of $\Sigma_2$ can be arbitrary.
\end{proof}

\begin{cor}\label{D+}
Suppose $\Sigma_D$ is the $D$-basis of some closure system. Consider $\Sigma_D^+\subseteq \Sigma_D$ obtained from $\Sigma_D$ by performing the following reductions:
\begin{itemize}
\item[(a)] Remove $A\rto x$, if $A\rto y$ and $y\rto x$ are also in $\Sigma_D$.
\item[(b)] Remove $z\rto x$, if $z\rto y$ and $y\rto x$ are also in $\Sigma_D$.
\end{itemize}
Let $\Sigma_1$ be a the proper ordering of binary part of $\Sigma_D^+$ given in Proposition \ref{cover}, and let $\Sigma_3$ be a subordering of this proper ordering on implications $y\rto x$ that appear in triples of $A\rto x, A\rto y, y\rto x$ of (a). Finally, let $\Sigma_2$ be some ordering of non-binary implications of $\Sigma_D^+$. 
Then $\sigma=\Sigma_1^\frown \Sigma_2^\frown\Sigma_3$ is the ordered direct sequence for the basis $\Sigma_D^+$. In particular, the length of this sequence is no longer than the length of the $D$-basis.
\end{cor}
\begin{proof}
Indeed, following the procedure of Lemma \ref{replacing} we can replace all $A\rto x$ from the triples $A\rto x, A\rto y, y\rto x$ in $\Sigma_D$ by the second copy of $y\rto x$ in additional binary part $\Sigma_3$ that follows the non-binary part of the $D$-basis.
\end{proof}

\begin{exm}
Given the $D$-basis of the closure system: $\Sigma_D=\langle 3\rto 2,2\rto 1, 3\rto 1, 45\rto 3,45\rto 2,45\rto 1\rangle$, we can produce a shorter basis $\Sigma_D^+=\{3\rto 2,2\rto 1,45\rto 3\}$ with the ordered direct sequence:
$\sigma=\langle 3\rto 2, 2\rto 1, 45\rto 3, 3\rto 2,2\rto 1\rangle$.
We note that $\Sigma_D^+$ is only half as long as $\Sigma_D$,, and its ordered direct sequence $\sigma$ has the same length as $D$-basis with optimized binary part but the size of $\sigma$ is smaller than that of the optimized $D$-basis.
\end{exm}


\section{Closure systems without $D$-cycles and the $E$-basis}\label{EB}

It turns out that the $D$-basis can be further reduced, when an additional property holds in a closure system $\langle S,\phi\rangle$. The results of this section follow closely the exposition given in \cite{FJN}, section 2.4.

We will write $xDy$, for $x, y \in S$, if $y \in Y$ for some minimal cover $Y$ of $x$. We note that the $D$-relation is a subset of the dependence relation $\delta$ from section~\ref{DRB}.

\begin{df}
A sequence $x_1,x_2,\dots,x_n$, where $n>1$, is called a $D$-\emph{cycle}, if $x_1Dx_2D\dots x_nDx_1$.
A finite closure system $\langle S,\phi \rangle$ is said to be \emph{without $D$-cycles} if it has no $D$-cycles. 
\end{df} 

We note that the lattices of closed sets of closure systems without $D$-cycles are known in lattice-theoretical literature as \emph{lower bounded}. 

For every $x \in S$, let $M(x)=\{Y\subseteq S: $Y$ \text{ is a minimal cover of } x\}$. The family $\phi(M(x))=\{\phi(Y): Y \in M(x)\}$ is ordered by set containment, so we can consider its minimal elements.
Let $M^*(x)=\{Y \in M(x): \phi(Y) \text{ is minimal in } \phi(M(x))\}$.

We will write $xEy$, for $x,y \in S$, if $y \in Y$ for some $Y \in M^*(x)$. According to the definition, if $xEy$ then $xDy$. On the other hand, the converse is not always true.

\begin{exm}\label{no-cycle}
Consider the closure system and its $D$-basis from Example \ref{Dbas}.
We note that this closure system has no $D$-cycles. We have three minimal covers of $6$: $15,24$ and $23$. Since $\phi(15)=S\setminus 2$, $\phi(24)=S$ and $\phi(23)=S\setminus 45$, we have only two of these covers in $M^*(6)$: $15$ and $23$.
Thus, while $6D4$, we do not have $6E4$.
\end{exm}

We now define two sequences of subsets of $S$, based on covers from $M(x)$ and $M^*(x)$, correspondingly.

Let $D_0=E_0=\{p \in S: p \in \phi(p_1,\dots, p_k) \text{ implies } p \in \phi(p_i) \text{ for some } i\leq k\}$.
If $D_k$ and $E_k$ are defined, then $D_{k+1}=D_k\cup \{s \in S: \text{ if } s\nor Y \text{ then } s\nor Z \text{ for some } Z \subseteq D_k, Z\ll Y\ \text { and } Z \in M(s)\}$. Similarly, 
$E_{k+1}=E_k\cup \{s \in S: \text{ if } s\nor Y \text{ then } s\nor Z \text{ for some } Z \subseteq E_k, Z\ll Y\ \text { and } Z \in M^*(s)\}$.
Apparently, $E_k\subseteq D_k$, for any $k$.
The following result is proved in \cite{FJN}, Theorem 2.51.

\begin{lem}\label{E=D} If $\langle S,\phi \rangle$ is a reduced closure system without $D$-cycles, then, for some $k$, $S=E_k=D_k$.
\end{lem}

As a consequence, we can often shorten the $D$-basis for a closure system without $D$-cycles.
We will say that $s\in S$ has $D$-rank $k=0$, if $s \in D_0$, and $k>0$, if $s \in D_k\setminus D_{k-1}$.
According to Lemma \ref{E=D}, every $s\in S$ in a closure system without $D$-cycles has a $D$-rank.

Recall that a basis is called \emph{aggregated} when all its premises are different. Every basis can be brought to the aggregated form by combining conclusions of all implications with the same premises.

\begin{thm}\label{Ebas}
Let $\langle S,\phi \rangle$ be a reduced closure system without $D$-cycles. Consider a subset $\Sigma_E$ of the $D$-basis that is the union of two sets of implications:
\begin{itemize}
\item[(1)] $ \{y \rightarrow x : x \in \phi(y)\}$,
\item[(2)] $\{X \rightarrow x: X \in M^*(x)\}$.
\end{itemize} 
Then 
\begin{itemize}
\item[(a)] $\Sigma_E$ is a basis for  $\langle S,\phi \rangle$. 
\item[(b)] $\Sigma_E$ is ordered direct.
\item[(c)] The aggregated form of $\Sigma_E$ is ordered direct.
\end{itemize}
\end{thm}
\begin{proof}
To begin with, it is not true that every cover of an element $x \in S$ 
refines to a cover in $M^*(x)$, so $\Sigma_E$ must be ordered more carefully
than $\Sigma_D$.  Nonetheless, mimicking the proof of Theorem~2.50 of \cite{FJN}, we can construct an order on $\Sigma_E$ that makes it an ordered direct basis.  This will be done for the aggregated $E$-basis, proving parts (a) and (c) simultaneously; part (b) then follows. 

Consider the aggregated form of $\Sigma_E$. Given an implication $X\rto Y$ in this basis, let $D^*(X\rto Y)$ be the maximal $D$-rank of elements in $X$, and $D_*(X\rto Y)$ be the minimal $D$-rank of elements in $Y$. Then $D^*(X\rto Y)<D_*(X \rto Y)$.
 
Order the implications following the rule: put the implications $x\rto Y$ first (aggregated form of binary part of $\Sigma_E$), and for the rest, if $D^*(X_1\rto Y_1)<D^*(X_2\rto Y_2)$
then $X_1\rto Y_1$ precedes $X_2\rto Y_2$ in the order. 

\begin{claim}
If $X_1\rto Y_1$ and $X_2\rto Y_2$ are in the aggregated $E$-basis, and $Y_1\cap X_2\not = \emptyset$, then $X_1\rto Y_1$ precedes $X_2\rto Y_2$.
\end{claim}

Indeed, take any $x \in Y_1\cap X_2$. If the $D$-rank of $x$ is $k$, then $D^*(X_2\rto Y_2)\geq k\geq D_*(X_1\rto Y_1)> D^*(X_1\rto Y_1)$. Hence, $X_1\rto Y_1$ will appear in the order before $X_2\rto Y_2$.

Now take any input set $Z$. We want to show that $\phi(Z)$ can be obtained when applying the aggregated basis in the described order. We argue by induction on the rank of an element $ z\in\phi(Z)\setminus Z$. 

If $z\in D_0$, then it only can be obtained via some implication $x\rto Y$, for some $x\in Z$, and $z \in Y$, and implications $x\rto Y$ form an initial segment in the ordered sequence of the basis. Now assume that it is already proved that all elements of $\phi(Z)\setminus Z$ of rank at most $k$ can be obtained in some initial segment of the sequence for the basis. If we have now element $z$ of rank $k+1$, then it can be obtained via an implication $X\rto Y$ with $X\subseteq \phi(Z)$, $z \in Y$, and $D^*(X)<k+1$. By the induction hypothesis, all elements in $X\subseteq \phi(Z)\setminus Z$ can be obtained via implications located in some initial segment of the sequence, and by the Claim above, all those implications precede $X\rto Y$. Thus, all implications producing elements of rank $k+1$  from $\phi(Z)$ will be located after the segment of the sequence producing all rank $k$ elements.   
\end{proof}

To illustrate the ordering of an $E$-basis, consider again the closure system given in Example \ref{Dbas}. As we know from Example \ref{no-cycle}, $\Sigma_E$ exists  and includes all implications of the $D$-basis, except $24\rto 6$.
Elements $1,2,4$ have $D$-rank $0$; elements $3,5$ have $D$-rank $1$, and $D$-rank of $6$ is 2. This allows to impose a proper ordering on implications of $\Sigma_E$ that turns it into ordered direct:\\
(1) $5\rto 4$, (2) $6\rto 3$, (3)  $3\rto 1$, (4)  $2\rto 1$, (5) $14\rto 3$, (6) $24\rto 5$, (7) $23\rto 6$, (8) $15\rto 6$. This basis is also aggregated.

\begin{prop}
Suppose $\Sigma_D=\{s_1,s_2,\dots,s_n\}$ is a $D$-basis of some reduced closure system $\langle S,\phi\rangle$  and $|S|=m$. It requires time $O(mn^2)$ to determine whether the closure system is without $D$-cycles, and if it is, to build its ordered direct basis $\Sigma_E$. 
\end{prop}
\begin{proof}
Since the $D$-relation is a subset of $S^2$, it will contain at most $m^2$ pairs. On the other hand, it is built from implications $X\rto x$, so the other upper bound for pairs in $D$-relation is $mn$. Evidently, the closure system is without $D$-cycles iff its $D$-relation can be extended to a linear order. There exists an algorithm that can decide whether $\langle S,D\rangle$ can be extended to a partial order on $S$ in time $O(m + |D|)$, see Theorem 11.1 in \cite{FJN}.
We will see below that the rest of the algorithm will take time $O(mn^2)$, which makes the total time also $O(mn^2)$.

Assuming the first part of algorithm provides a positive answer and there are no $D$-cycles, we proceed by finding the ranks of all elements.
It will take at most $n$ operations to find set $D_0$: include $p$ into $D_0$, if it does not appear as a conclusion in any (non-binary) implication $X\rto x$ of the $D$-basis, where $x\nor X$. If the system is without $D$-cycles, then $\pi_\Sigma(D_0)\setminus D_0$ gives elements of rank 1, $\pi_\Sigma^2(D_0)\setminus \pi_\Sigma(D_0)$ elements of rank 2, etc. Note that $\pi_\Sigma(X)$ is defined in the beginning of section 6. 
Computation of $\pi_\Sigma(X)$ requires time $O(mn)$, since $\Sigma=\Sigma_D$ in our case has $n$ implications, and checking that the premise of each implication is a subset of $X$ takes time $O(m)$. After at most $m$ iterations of $\pi$ on $D_0$, one would obtain the whole $S$, whence, $O(m^2n)$ operations are needed to obtain the ranks of all elements from $S$. Assuming that $m\leq n$ in most closure systems, this time will not beat $O(mn^2)$.

 It remains to decide which implications from the $D$-basis should remain in the $E$-basis. To that end, for each element $x \in S$ we need to compare the closures $\phi(X)$ of subsets $X$, for which $X \rto x$ is in the $D$-basis, and choose for the $E$-basis those that are minimal.
There is at most $n$ implications $X\rto x$, for a given $x \in S$, and the closure $\phi(X)$, for each such $X$, can be found in $O(s(\Sigma_D))$ steps. It will take time $O(n^2)$ to determine all minimal subsets among $O(n)$ given subsets $\phi(X)$, associated with fixed $x \in S$. Hence, it will require time $O(mn^2)$ for all $x \in S$. 

The size of the $E$-basis will be at most $n$, and it will take time $O(n^2)$ to order it with respect to the rank of elements, per Theorem \ref{Ebas}.
\end{proof}


When a closure system has $D$-cycles, the subset $\Sigma_E$ of $\Sigma_D$, defined in Corollary \ref{Ebas}, may not form a basis. 

\begin{exm}\label{ex22}
Consider $S=\{1,2,3,4\}$ and a closure operator defined by the $D$-basis
\[
13\rto 2, 24\rto 3,14\rto 2,14\rto 3.
\]
This closure system has the cycle $2D3D2$. It is easy to verify that $\Sigma_E$ has only $13\rto 2$ and $24\rto 3$, so the last two implications from the $D$-basis cannot be recovered from $\Sigma_E$.
\end{exm}

\begin{figure}[htbp]
\begin{center}
\includegraphics[height=2.2in,width=6.0in]{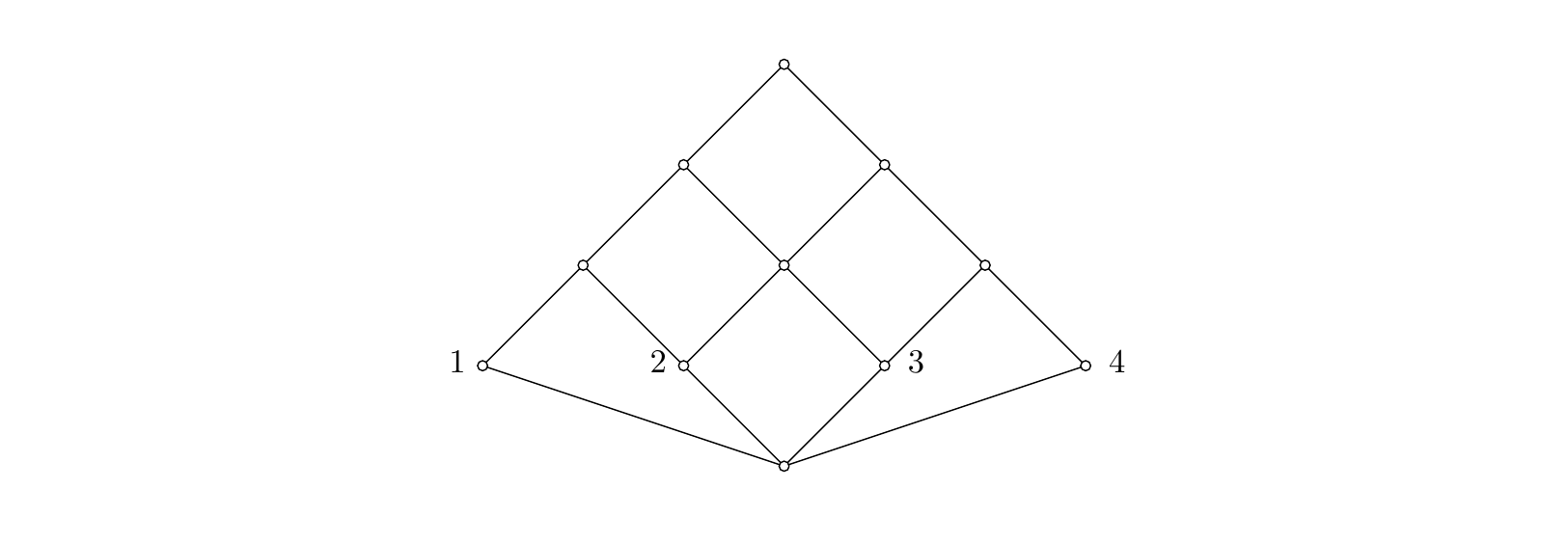}
\caption{Example \ref{ex22}}
\end{center}
\end{figure}

Further results about closure systems without D-cycles, and more generally systems whose closure lattice is join semidistributive, will be presented in~\cite{JK2}.

\section{$D$-basis versus Duquenne-Guigues canonical basis}\label{D and DG}

We recall the definition of the canonical basis introduced by V.~Duquenne and J.L.~Guigues in \cite{DG}, see also \cite{CM03}.  This applies to arbitrary closure systems, not just reduced ones.

\begin{df}\label{canon}
The canonical basis of a closure system $(S,\phi)$ consists of implications $X\rightarrow Y$ for $X,Y \subseteq S$, that satisfy the following properties:
\begin{itemize}
\item[(1)] $X \subset \phi(X)=Y$;
\item[(2)] for any $\phi$-closed set $Z$, either $X \subseteq Z$ or $Z \cap X$ is $\phi$-closed;
\item[(3)] if $W\subseteq X$, $\phi(W)=Y$ and $W$ satisfies (2) in place of $X$, then $W=X$. 
\end{itemize}
\end{df}

The subsets $X\subset S$ with properties (1) and (2) are usually called \emph{quasi-closed}, see \cite{CM03}.
The meaning of (2) is that adding $X$ to the family of closed sets of $\phi$ produces the family of closed sets of another closure operator. Property (3) indicates that among all quasi-closed subsets with the same closure one needs to choose the minimal ones.
This basis is called \emph{canonical}, since it is minimal, in that no implication can be removed from it without altering $\phi$, and every other minimal implicational basis for $\phi$ can be obtained from it. In particular, no other basis can have a smaller number of implications. Note that here the implications are of the form $X\rto Y$, where $Y$ is not necessarily a one-element set. We will also call it the \emph{D-G basis}, to distinguish from canonical unit direct basis. 
 
To bring this basis in comparison with other bases discussed in this paper, each implication $X \rightarrow Y$ may be replaced by set of implications $X\rightarrow y$, $y \in Y\setminus X$. We will call this modification of the canonical basis the \emph{unit D-G basis}.

In many cases the canonical basis may be turned into an ordered direct basis. 

\begin{exm}
Consider again the closure system from Example \ref{Dbas}. The canonical basis is
\[
2\rto 1,3\rto 1,5\rto4,6\rto 3,6\rto 1, 14\rto 3, 123\rto 6, 1345\rto 6, 12346\rto 5.
\]
Besides, it is ordered direct in the given order.
\end{exm}

In general, though, the canonical basis cannot be ordered so that it becomes direct. Thus, it is not ordered direct. The following two examples were uncovered by running a computer program and checking about a million of various closure systems on $5$- and $6$-element sets.
The first example demonstrates a closure system, where the canonical basis cannot be ordered, while the unit expansion of this basis does admit an ordering to make it direct. The second example shows that some canonical bases cannot be ordered in either form. 

\begin{figure}[htbp]
\begin{center}
\includegraphics[height=2.2in,width=6.0in]{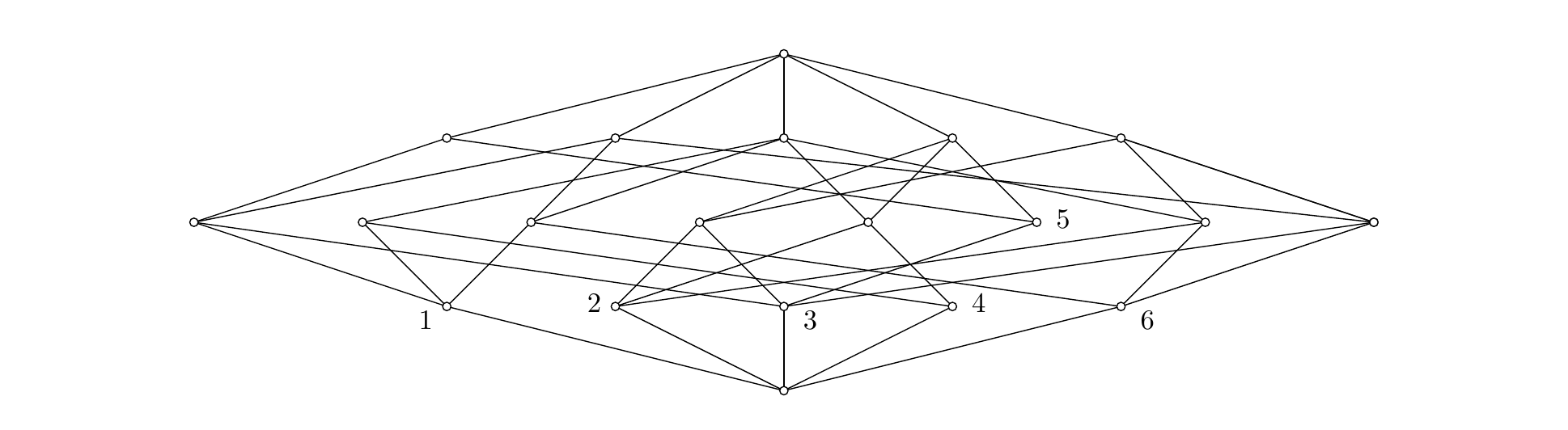}
\caption{Example \ref{ex66}\label{pic66}}
\end{center}
\end{figure}

\begin{exm}\label{ex66}
\end{exm}
Let $\langle S,\phi\rangle$ be a closure system on $S=\{1,2,3,4,5,6\}$, given by the family of closed sets:
$\{\emptyset, 1,2,3,4,6, 3 6, 2 6, 1 3, 2 4,  1 4, 3 5, 2 3, 1 6, 1 3 5,  1 3 6, 2 36, 1 2 4 6, 2 3 4 5, S\}$.
The lattice representation of this system is given in Figure \ref{pic66}.

Then the canonical basis is $5\rto 3, 34\rto 25, 12\rto 46, 46\rto 12, 235\rto 4, 356\rto 124$.
It is easy to show that this basis cannot be ordered. Indeed, in order to obtain $\phi(145)=S$ in one application of canonical basis, one would need to put $5\rto 3$ first, then $34\rto 25$, followed by $12\rto 46$.
On the other hand, $\phi(123)=S$, too, and the only implication applicable to $123$ is $12\rto 46$, but it comes after $34\rto 25$, and one cannot obtain $5$ in the closure otherwise.

As was mentioned, the unit expansion of this canonical basis is still ordered direct: one would need to place implications $12\rto 4$ and $12\rto 6$ around $34\rto 2$ and $34\rto 5$, thusly:
$5\rto 3, 12\rto 4, 34\rto 2, 34\rto 5, 12\rto 6, 46\rto 2, 46\rto 1, 235\rto 4, 356\rto 1, 356\rto 2, 356\rto 4$.

As always, the $D$-basis is ordered direct in both forms: in its original unit form, and in the aggregated form. 
For example, the aggregated form of $D$-basis in this example is
$5\rto 3, 34\rto 25, 12\rto 46, 46\rto 12, 25\rto 4, 56\rto 124, 123\rto 5, 134\rto 6$.

One needs to run the canonical basis two times to ensure the closure of arbitrary subset, i.e., apply $6 \cdot 2=12$ implications, versus only $8$ implications of the aggregated $D$-basis. 
In the unit form, the canonical basis has $11$ implications and the $D$-basis has $13$, but the ordering of the canonical basis requires special care.

\begin{rem} \label{rem29}
\end{rem}  Example \ref{ex66} also shows that the $D$-basis, unlike the canonical basis, can be \emph{redundant} (even in its aggregated form): this means that some implications can be removed, and the remaining ones still define the same closure system. In the $D$-basis of our example, both implications $123\rto 5, 134\rto 6$ can be removed, since they follow from $34\rto 25, 12\rto 46$. On the other hand, the basis without these two implications is no longer ordered direct.

The following example shows that the canonical basis might be un-orderable in either form.

\begin{figure}[htbp] 
\begin{center}
\includegraphics[height=2.2in,width=6.0in]{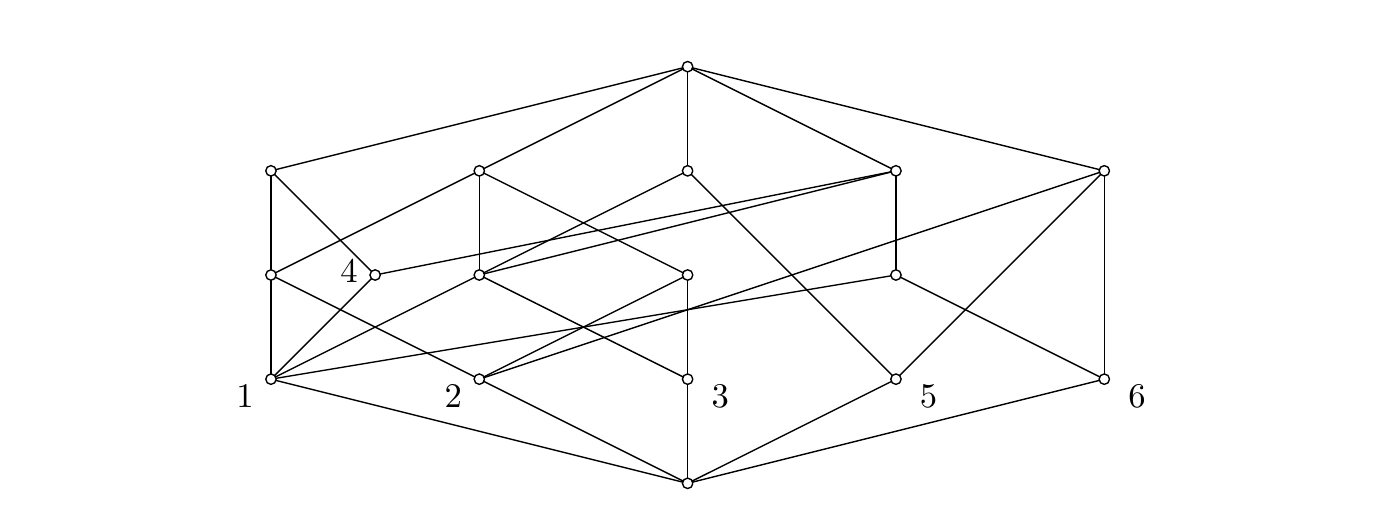}
\caption{Example \ref{Example 67}}\label{pic67}
\end{center}
\end{figure}

\begin{exm}\label{Example 67}
\end{exm}
Let $\langle S,\phi\rangle$ be a closure system on $S=\{1,2,3,4,5,6\}$, given by the family of closed sets:
$\{\emptyset, 1,2,3,5,6, 12, 13, 14, 16, 23, 123, 124, 135, 256, 1346, S\}$.
The lattice representation of this system is given in Figure \ref{pic67}.

The canonical basis has 9 implications:\\
$4\rto 1, 15\rto 3, 35\rto 1, 25\rto 6, 56\rto 2, 26 \rto 5, 36\rto 14, 134\rto 6, 146\rto 3$.\\
There is a single implication $36\rto 14$ that can be expanded to two unit implications $36\rto 1$ and $36\rto 4$.

The proof that the unit expansion of canonical basis cannot be ordered to make it direct, follows from consideration of the next three closures:
\begin{itemize}
\item $45\rto 145\rto 1345\rto 13456 \rto S$, hence, $134 \rto 6$ should be placed later than $15\rto 3$.
\item $1234 \rto 12346 \rto S$, hence $26\rto 5$ should be placed later than $134 \rto 6$.
\item $126 \rto 1256 \rto 12356 \rto S$, hence $15\rto 3$ should be placed later than $26\rto 5$, which contradicts the combination of the previous two items.
\end{itemize}

For comparison, the aggregated $D$-basis has 15 implications:\\
$4\rto 1, 45\rto 26, 36\rto 14, 34\rto 6, 15\rto 3, 46 \rto 3, 35\rto 1, 25\rto 6, 26\rto 5, 56\rto 2, 126 \rto 34, 235 \rto 4, 156 \rto 4, 234 \rto 5, 125 \rto 4$.

Thus, one run of the aggregated $D$-basis (15 implications) wins over two runs (18 implications) of the canonical basis.
In unit expansions: $D$-basis (18 implications) still wins over two runs (20) of canonical basis. 

In this example, the $D$-basis is $4$ implications shorter than the canonical unit direct basis, which has $22$ implications.\\

Our earlier analysis of the binomial part of the $D$-basis in Proposition \ref{cover} carries over to a partial optimization of the canonical basis.

\begin{prop}\label{binary} The binary part of the unit expansion of the D-G canonical basis of any reduced closure system coincides with the binary part of the $D$-basis (or, $E$-basis, if it exists) of the same system.
\end{prop}
\begin{proof}
We recall that the binary part of the $D$-basis of closure system $\langle S,\phi\rangle$ consists of implications $y\rto x$, where $x \in \phi(y)\setminus y$. This implies that $\{y\}$ is not a $\phi$-closed set. Besides, it is a quasi-closed set, since the intersection of $\{y\}$ with any $\phi$-closed set is either $\{y\}$ or $\emptyset$. Evidently, $\{y\}$ will be the minimum quasi-closed set with the closure $\phi(y)$. Hence, $\{y\}\rto \phi(y)\setminus y$ should be an implication in the canonical basis. Evidently, the unit expansion of  $\{y\}\rto \phi(y)\setminus y$ gives all the implications in the $D$-basis with the premise $y$. Vice versa, every implication in the canonical basis of the form $y \rto Y$ implies that $Y=\phi(y)\setminus y$.  Hence, $y\rto y'$ for $y'\in Y$ should appear in the $D$-basis.
\end{proof}

The following statement is an immediate consequence of Proposition \ref{cover} and Proposition \ref{binary}.
We recall that $L$ stands for the lattice of closed sets of $\langle S,\phi\rangle$, and $(J(L),\leq)$ is a partially ordered set of join-irreducible elements of $L$.

\begin{cor}
Let $\Sigma_C$ be the canonical basis of $\langle S,\phi\rangle$, where $|S|=m$. Let $\Sigma_C^b\subseteq \Sigma_C$ be a binary part of $\Sigma_C$, and let $n$ be the number of implications in the unit expansion of $\Sigma_C^b$. Then an algorithm that requires $O(mn+n^2)$ time will replace each implication $y\rto Y$ in $\Sigma_C^b$ by $y\rto Y'$, $Y'\subseteq Y$, where $\phi(y)$ covers $\phi(y')$ in $(J(L), \leq)$, for each $y' \in Y'$.
Let $\Sigma_C'$ be this new set of implications. The algorithm will also put an appropriate order on $\Sigma_C'$ in such a way that $\rho_{\Sigma_C}=\rho_{\Sigma_C'}$.
\end{cor}

Thus, the optimization of the canonical basis inspired by Proposition \ref{cover} is in the form of a possible size reduction of some implications.\\

We finish this section with a comparison of the canonical D-G basis with the $D$-basis on some illustrative examples.
We consider one particular type of closure systems for which the description of the canonical basis is easy.
The closure system $\<S,\phi\>$ is called \emph{a convex geometry}, if $\phi$ satisfies the anti-exchange axiom: if $x \in \phi(C\cup\{y\})$ and $x \notin C$, then $y \notin \phi(C\cup\{x\})$, for all $x \ne y$ in $S$ and all closed $C \subseteq S$.

For any closed set $X$ in a convex geometry, the set of extreme points of $X$ is defined as $Ex(X)=\{x \in X: x \not \in \phi(X\setminus x)\}$. It is well-known that, in every convex geometry, $X=\phi(Ex(X))$. The equivalent statement in the framework of lattice theory is that every element $Y$ in the closure lattice of a finite convex geometry has unique representation as a join of join irreducible elements: $Y=\bigvee Y_i$, so that none of $Y_i$ can be removed (such representation is called \emph{irredundant join decomposition} of $Y$); see, for example, \cite{AGT}.

An important example of convex geometry is $Co(R^n,A)$, where $A$ is a finite set of points in $R^n$, and $Co(R^n,A)$ stands for geometry of convex sets relative to $A$. In other words, the base set of such closure system is $A$, and closed sets are subsets $X$ of $A$ with the property that whenever point $a \in A$ is in convex hall of some points from $X$, then $a$ must be in $X$ (see more details of the definition, for example, in \cite{AGT}).

\begin{lem}\label{quasiclosed}
If $Y$ is the premise of an implication from the canonical basis of some convex geometry $Co(R^n,A)$, then
$Y$ is the set of extreme points of a closed set $\phi(Y)$ such that every subset of $Y$ is closed.
\end{lem}
\begin{proof}
Evidently, the premise $Y$ of every implication of the canonical basis contains $ex(\phi(Y))$. Moreover, 
$Co(R^n,A)$ satisfies the $n$-Carath\'eodory property, see \cite{A08}, which means that $|Ex(\phi(Y))|\leq n+1$.
Suppose there exists an implication $Y\rto z$ in canonical basis with $Ex(Y)=\{y_1,y_2,\dots, y_{n+1}\}$, $z \in \phi(Y)$ and $z \not \in \phi(Y')$ for every $Y'\subset ex(Y)$. We claim that every $Y_i=Y\setminus \{y_i\}$ is closed. Indeed, suppose w.l.o.g. that $x \in \phi(Y_{n+1})=\phi(y_1,\dots, y_n)$, $x \not \in \{y_1,\dots, y_n\}$. Then simplex generated by $y_1,\dots, y_{n+1}$ is split into simplices generated by $X_1=\{x,y_2,\dots, y_{n+1}\}$, $X_2=\{y_1,x,y_3,\dots, y_{n+1}\}$,\dots, $X_n=\{y_1,y_2,\dots, y_{n-1},x,y_{n+1}\}$. Then $z$ must be in one of those simplices, say, $z \in \phi(X_i)$. Since $Y$ is quasi-closed, $\phi(Y_{n+1}) \subset \phi(Y)$ implies $x \in Y$, and $\phi(X_i)\subset \phi(Y)$ implies $z \in Y$, a contradiction with the assumption $z \not \in Y$.

Similar argument applies for any other $Y$ with $E
x(\phi(Y)) <n+1$.  
\end{proof}

In the next two examples, we consider convex geometries of the form $Co(R^2,A)$ and compare the canonical bases and $D$-bases.

\begin{exm}
If $A$ is a set of points in general position, i.e., no three points are on a line, then the $D$-basis and canonical basis of convex geometry $Co(R^2,A)$ are the same.
\end{exm}

Due to the $3$-Carath\'eodory property, all covers can be reduced to covers by three elements. So the $D$-basis consists of implications $abc\rightarrow x$, for all triangles $abc$ that have $x$ inside.

Now, for any relatively convex subset $X \subseteq A$, $Ex(X)$ consists of the vertices of a convex polygon that holds all the points of $X$ inside. If there exists $y \in X\setminus Ex(X)$, then there are $a,b,c \in X$ with $y \in \phi(a,b,c)$. Hence, $\{a,b,c\} \subseteq X$ is not closed. This will not contradict Lemma \ref{quasiclosed}, only if $X=\{a,b,c\}$.
Thus, the only implications in the canonical basis are $abc \rightarrow x$, where $x$ is inside triangle $abc$. 
It follows that the $D$-basis and canonical basis are the same.

\begin{exm}\label{DnotCanonic}
If $A$ is a set of points that is not in general position, then the canonical basis of $Co(R^2,A)$ is a proper subset of the $D$-basis.
\end{exm}

Indeed, consider a point configuration of 5 points: $a,b,c$ form a triangle, $x$ is inside the triangle, and $d$ is on the side $ab$, so that $x$ is also inside triangle $dbc$.
The $D$-basis is $ab\rightarrow d$, $abc\rightarrow x$, $bcd \rightarrow x$, while the canonical basis is $ab \rightarrow d$, $bcd\rightarrow x$.

Note that $abc$ cannot be a premise of an implication in the canonical basis due to Lemma \ref{canon}, since the subset $ab$ is not closed. \\

We note that Lemma \ref{quasiclosed} is not true for arbitrary convex geometries. 
\begin{exm}
\end{exm}
Take convex geometry  $(\{a,b,c,d,x\},\phi)$ of Example \ref{DnotCanonic}. Adding another closed set $\{b,c,d\}$ will result in a new convex geometry $(\{a,b,c,d,x\},\psi)$ with the canonical basis $ab\rto d, abcd\rto x$. Note that in implication $abcd\rto x$, $d \not \in Ex(abcd)$, and subset $ab$ is not closed.

\section{Testing the performance of $D$-basis}\label{performance}

The performance of $D$-basis in comparison with the D-G unit basis and canonical unit direct basis was tested on $300,000$ randomly generated closure systems on base sets of $6$ and $7$ elements. The closed sets in these systems were generated by taking $3$ to $8$ arbitrary subsets of the domain, the intersection of all combinations of these sets, the empty set, and the domain itself.

The computation of the closure of random input set $X$ was implemented, for the D-G unit basis, according to the \emph{folklore} algorithm, which essentially makes the computation of $\pi(X)$, $\pi^2(X)$, etc., on its consecutive loops.
This algorithm is presented, for example, as Algorithm $0$ in section 2 of \cite{W95}, also see our discussion of this algorithm in comparison with the forward chaining algorithm in section \ref{forward}.  
Based on this algorithm, computing the closure of an input set using the D-G basis will always take at least two passes: the final pass produces nothing and exists solely to determine that the ability of the basis to expand the given set has been exhausted.

In contrast, the computation of the closure of any input set, by the $D$-basis or canonical unit direct basis, is done simply in one loop of such algorithm. The data collected reflects the number of implications attended in the run of each algorithm. Thus, with respect to these two bases, it makes a comparison of their length.

\begin{table}
\label{table6}
\begin{tabular}{|l|l|l|l|l|}
\hline
Closed Sets & D-G Unit Basis & Direct Optimal Basis & D-Basis \\ \hline
5 & 32.37 & 18.71 & 16.45 \\ \hline
10 & 22.59 & 16.24 & 12.23 \\ \hline
15 & 21.22 & 15.54 & 12.47 \\ \hline
20 & 18.13 & 13.06 & 11.45 \\ \hline
25 & 15.54 & 11.09 & 10.34 \\ \hline
30 & 11.70 & 7.96 & 7.65 \\ \hline
\end{tabular}
\caption{Average implications checked to expand an arbitrary 3-element set in a length 6 domain}
\end{table}

In the testing on domain length $6$, with inputs sets of length $3$, the D-G unit basis cycled through, on average, $22.9$ implications before returning the closure. By comparison, the direct canonical (optimal) basis took $15.8$ such steps and the $D$-basis took only $12.7$ checks on average. Due to their ordered directness, the number of implications checked in the direct optimal and $D$-basis was equivalent to the number of implications they contained. 

\begin{figure}[htbp]
\label{graph6}
\begin{center}
\includegraphics[width=\textwidth]{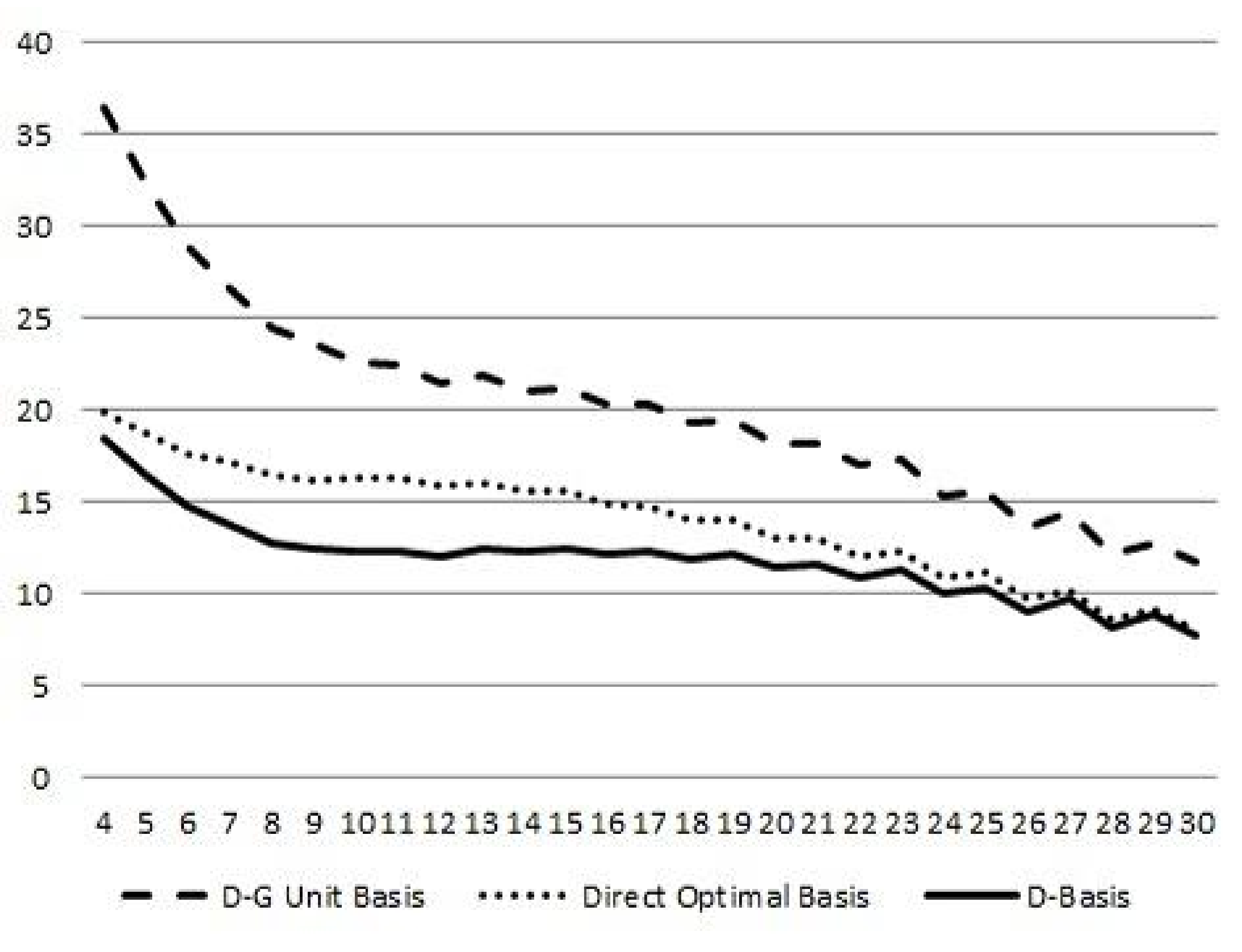}
\end{center}
\caption{Implications checked (Y-Axis) by the number of closed sets in the basis (X-axis) for a domain of size 6.}
\label{domain6graph}
\end{figure}

It was observed that the efficiency gap between the direct and indirect bases was greatest when there were fewer closed sets, meaning that more subsets could be expanded through the bases' implications. This relation is shown in Figure \ref{domain6graph}. 

We saw similar results on bases of domain length $7$. There, we once again saw the convergence of the $D$-Basis and direct optimal as the number of closed sets approached either extreme, with a more pronounced gap in between.

\begin{table}
\label{table7}
\begin{tabular}{|l|l|l|l|l|}
\hline
Closed Sets & D-G Unit Basis & Direct Optimal Basis & D-Basis  \\ \hline
5 & 46.73 & 27.70 & 23.57  \\ \hline
10 & 33.74 & 26.26 & 17.92 \\ \hline
15 & 32.11 & 26.80 & 18.59 \\ \hline
20 & 31.01 & 25.68 & 19.43 \\ \hline
25 & 29.77 & 23.99 & 19.66 \\ \hline
30 & 26.71 & 20.64 & 17.73 \\ \hline
\end{tabular}
\caption{Average implications checked to expand an arbitrary 3-element set in a length 7 domain.}
\label{table7}
\end{table}


There were $33.8$ checks on average for the D-G unit basis, and $26.0$ and $19.0$ for the direct optimal and $D$-basis, respectively, see Table \ref{table7}.

\begin{figure}[htbp]
\label{graph7}
\begin{center}
\includegraphics[width=\textwidth]{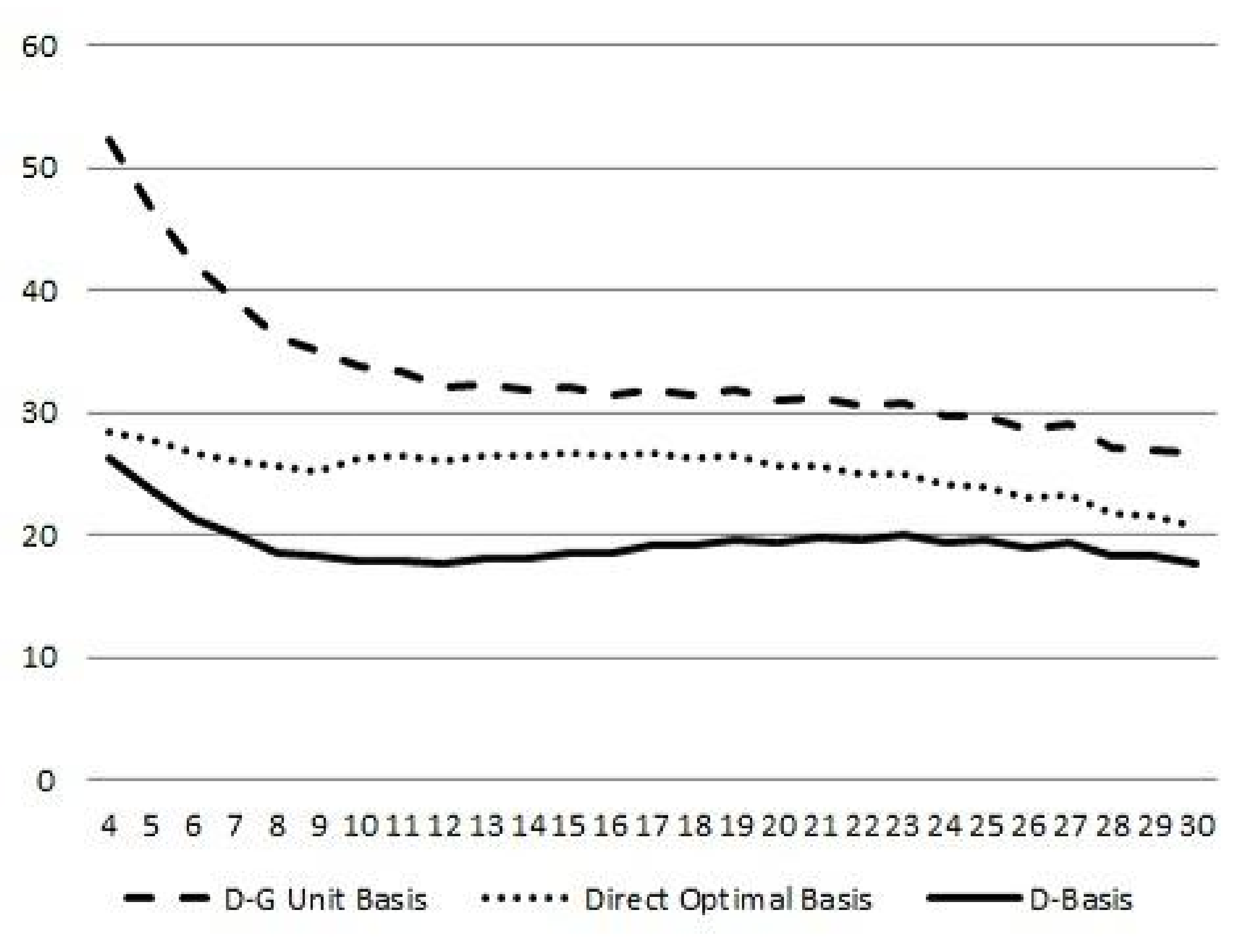}
\end{center}
\caption{Implications checked (Y-Axis) by the number of closed sets in the basis (X-axis) for a domain of size 7.}
\end{figure}

\emph{Acknowledgments.} We are grateful to several people who read the draft of this paper at different stages of its preparation and made valuable comments: Gy\"orgy Tur\'an, Vincent Duquenne, Marina Langlois, Ralph Freese, Marcel Wild and Karell Bertet, also students of Yeshiva College Joshua Blumenkopf and Jeremy Jaffe. 
Major advancements of this paper were done during the first author's visits to University of Hawai'i in 2010-2011, supported by the AWM-NSF Mentor Travel grant. The welcoming atmosphere at the Department of Mathematics at UofH is greatly appreciated. We thank our anonymous referees for making the numerous comments, pointing to a few references and fixing some of our omissions.

\end{document}